\documentclass[fside]{article}
\usepackage[utf8]{inputenc}
\usepackage[english]{babel}
\usepackage{lipsum}
\usepackage{bm}
\usepackage{upgreek}
\usepackage{graphicx}
\usepackage{stmaryrd}
\usepackage[linesnumbered, ruled]{algorithm2e}

\SetKwRepeat{Do}{do}{while}%
\usepackage{listings}
\usepackage{graphicx}
\usepackage{adjustbox}
\usepackage{etoolbox}
 \usepackage{multirow}
\usepackage[numbers,sort&compress]{natbib}  

\usepackage[framemethod=TikZ]{mdframed}  
\usepackage{url}  
\usepackage{subcaption}  
\usepackage{array}
\usepackage{booktabs}
\usepackage{makecell}
\usepackage{amsmath}
\usepackage{float}
\usepackage{titlesec}
\newcolumntype{C}{>{\centering\arraybackslash}X}
\newcolumntype{R}{>{\raggedleft\arraybackslash}X}
\newcolumntype{L}{>{\raggedright\arraybackslash}X}

\usepackage{changepage}
\usepackage{authblk}
\usepackage{threeparttable}

\usepackage{xcolor}  
\usepackage{geometry}
\usepackage{algpseudocode}

\usepackage{amsmath}

\usepackage{amsthm,amsmath,amssymb}
\usepackage{mathrsfs}
\usepackage{microtype} 

\definecolor{MyColor1}{rgb}{0,0,0}

\usepackage{titlesec}
\usepackage{sectsty}
\sectionfont{\color{MyColor1}}  
\subsectionfont{\color{MyColor1}}  

\usepackage{caption}
\usepackage{subcaption}

\usepackage{tocloft}
\makeatletter
\let\reftagform@=\tagform@
\def\tagform@#1{\maketag@@@{(\ignorespaces\textcolor{red}{#1}\unskip\@@italiccorr)}}
\renewcommand{\eqref}[1]{\textup{\reftagform@{\ref{#1}}}}
\makeatother

\usepackage{setspace} 

\usepackage[utf8]{inputenc}  

\title{Traffic Flow Data Completion and Anomaly Diagnosis via Sparse and Low-Rank Tensor Optimization}
\author{Junxi Man$^*$, Yumin Lin\thanks{First author and second  author contribute equally to this work.} ~and Xiaoyu Li\thanks{Corresponding author. Email: 21118027@bjtu.edu.cn. Her research work was supported by the Fundamental Research Funds for the Central Universities (2023YJS073) and the National Natural Science Foundation of China (12271022).}}
\affil{School of Mathematics and Statistics, Beijing Jiaotong University}

\date{} 
\usepackage{setspace} 
\begin{document}
\setlength{\baselineskip}{1\baselineskip}
\maketitle

\renewcommand{\abstractname}{%
    \fontsize{13}{16}\selectfont 
    \textbf{Abstract} 
}
\begin{abstract}
\fontsize{12}{15}\selectfont
\setstretch{1.5}
Spatiotemporal traffic time series, such as traffic speed data, collected from sensing systems are often incomplete, with considerable corruption and large amounts of missing values. A vast amount of data conceals implicit data structures, which poses significant challenges for data recovery issues, such as mining the potential spatio-temporal correlations of data and identifying abnormal data. In this paper, we propose a Tucker decomposition-based sparse low-rank high-order tensor optimization model (TSLTO) for data imputation and anomaly diagnosis. We decompose the traffic tensor data into low-rank and sparse tensors, and establish a sparse low-rank high-order tensor optimization model based on Tucker decomposition. By utilizing tools of non-smooth analysis for tensor functions, we explore the optimality conditions of the proposed tensor optimization model and design an ADMM optimization algorithm for solving the model. Finally, numerical experiments are conducted on both synthetic data and a real-world dataset: the urban traffic speed dataset of Guangzhou. 
Numerical comparisons with several representative existing algorithms demonstrate that our proposed approach achieves higher accuracy and efficiency in traffic flow data recovery and anomaly diagnosis tasks.
\end{abstract}



\section{INTRODUCTION}

With the widespread acquisition of large-scale traffic data and the extensive application of Intelligent Transportation Systems (ITS), effective traffic data completion and anomaly diagnosis have become critical challenges in traffic data research. Traditional traffic data completion methods mainly rely on the global low-rankness assumption, recovering missing data by leveraging matrix or tensor decomposition techniques. For instance, Liu et al.\cite{halrtc} proposed the HaLRTC model, enhancing completion accuracy by minimizing the sum of nuclear norms (SNN) across all modes. Chen et al.\cite{LRTCTNN} introduced the truncated nuclear norm (TruNN) and developed the LRTC-TruNN model, further optimizing low-rank approximation. However, these methods often overlook the local consistency of traffic data, particularly spatial and temporal correlations, leading to suboptimal completion performance for complex traffic network data and insufficient accuracy and robustness in anomaly detection tasks. To address this, Chen et al. \cite{LATC} incorporated autoregressive models to characterize temporal continuity while preserving global low-rankness.

To better capture local features in traffic data, recent studies have explored integrating regularization terms into traditional matrix or tensor decomposition frameworks. For example, Yu et al. \cite{trmf} proposed the Temporal Regularized Matrix Factorization (TRMF) model, which imposes graph regularization on factor matrices to handle missing values in time-series forecasting. Wang et al. \cite{omega} introduced the OrTC model based on CP decomposition, utilizing $\ell_{2,1}$-norm to smooth slice-wise anomalies and $\ell_2^2$-norm  to regularize  noise. Lyu et al. \cite{toep} developed a Tucker decomposition-based model with $\ell_1$-norm constraints on the core tensor and error tensor.
\begin{figure}[H]
\centering
{\includegraphics[width=.3\textwidth]{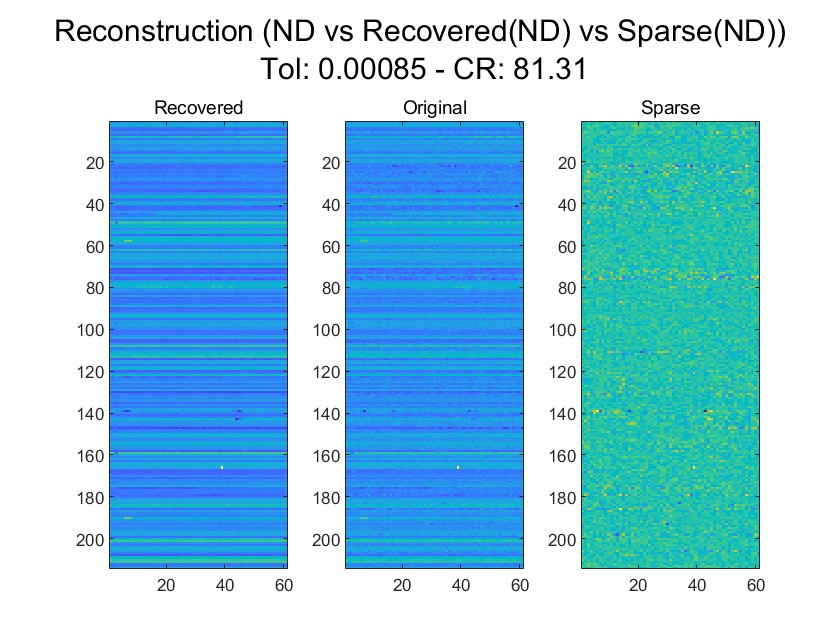}}
{\includegraphics[width=.26\textwidth]{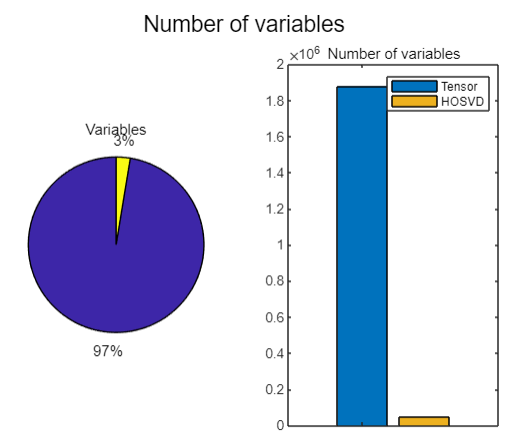}}
\caption{Application of Tucker decomposition in real traffic dataset. On the left side, from left to right: raw data, low-rank structure, sparse structure; on the right side: compression ratio.}\label{fig:4}
\end{figure}

\begin{figure}[H]
    \centering
    \includegraphics[width=0.45\linewidth]{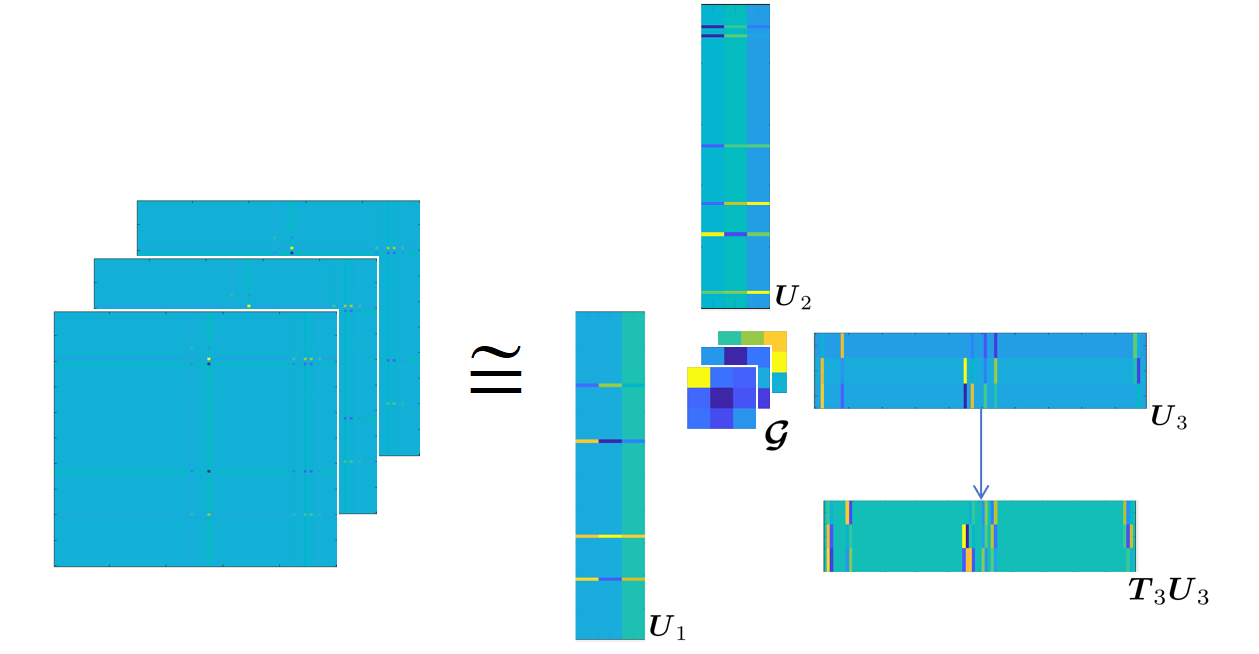}
    \caption{Core tensor and factor matrices generated by Tucker decomposition}
    \label{fig:nya}
\end{figure}

In this work, we structure traffic data into a third-order tensor $\mathcal{X} \in \mathbb{R}^{I \times J \times K}$ (sensors $\times$ time interval $\times$ day) and aim to decompose it into low-rank and sparse tensor components. Since the traffic tensor data has potentially complex spatiotemporal correlation, a powerful tool is to use tensor decomposition to mine spatiotemporal correlation, such as Tucker decomposition. As illustrated in Figure~\ref{fig:4}, Tucker decomposition effectively captures the global low-rankness inherent in real-world traffic data. Considering the periodicity and propagation of traffic spatio-temporal data, a strategy is adopted to enhance temporal and spatial similarity with low-rank component being encoded via Toeplitz matrices. Figure~\ref{fig:nya} depicts a synthetic data set with global low-rankness and local continuity, together with its Tucker decomposition components: three factor matrices $\boldsymbol{U}_1$, $\boldsymbol{U}_2$, $\boldsymbol{U}_3$, core tensor $\mathcal{G}$, and the transformed matrix $\boldsymbol{T}_3\boldsymbol{U}_3$ obtained by multiplying $\boldsymbol{U}_3$ with a Toeplitz matrix $\boldsymbol{T}_3$. Notably,  $\boldsymbol{T}_3\boldsymbol{U}_3$ exhibits numerous all-zero rows, which inspires us to formulate a regularization term that minimizes the number of non-zero rows in the Toeplitz-transformed factor matrices, thereby enforcing local continuity in the low-rank component of real-world data.

For the anomaly component, we leverage the $\ell_0$-norm to characterize the sparsity of the anomaly tensor $\mathcal{R}$, motivated by the sporadic nature of anomalous events. Furthermore, since anomalies often propagate within localized spatial-temporal neighborhoods (e.g., traffic congestion spreading to adjacent regions and time slots), the anomaly tensor exhibits block-wise sparsity. To model this property, we apply Toeplitz matrices $\boldsymbol{T}_l$ and $\boldsymbol{T}_r$ to the mode-1 unfolding of $\mathcal{R}$, denoted as $\boldsymbol{T}_l \mathcal{R}_{[1]} \boldsymbol{T}_r^\top$, thereby encoding short-term spatial-temporal propagation patterns of anomalies.

We propose a \textbf{T}ucker decomposition-based \textbf{S}parse \textbf{L}ow-Rank high-Order \textbf{T}ensor \textbf{O}ptimization model (TSLTO) framework, which systematically addresses high-order tensor analysis through three integrated mechanisms: (1) Global low-rank structure characterization via Tucker decomposition, (2) Anomaly-aware sparse regularization through $\ell_0$-norm regularized optimization, and (3) Spatiotemporal continuity preservation employing Toeplitz-constrained tensor factorization to maintain spatial-temporal consistency in both low-rank representations and sparse residuals. This unified framework not only effectively models the low-rank structure of regular traffic patterns but also precisely identifies sparse anomalies while maintaining local continuity, leading to enhanced missing data recovery and anomaly detection accuracy. To solve this non-convex optimization problem, we adopt the Alternating Direction Method of Multipliers (ADMM) algorithm, which balances sparsity and continuity constraints efficiently, ensuring computational scalability without compromising precision.

Our contributions are threefold: 
(i) A novel tensor completion sparse optimzation model  integrating Tucker low-rankness with orthogonal constraint for capturing the global spatiotemporal correlations, $\ell_0$-norm for describing the sparsity of the anomalies data, and Toeplitz-based continuity for mining the local spatiotemporal correlations;  
(ii) An efficient ADMM algorithm based on the a curvilinear search method with BB steps  is designed to solve the  non-convex and non-smooth  optimization model;  
(iii) Empirical and theoretical insights into leveraging spatial-temporal correlations for traffic data imputation and anomaly diagnosis.  


The remaining parts of this paper are organized as follows. Section 2 includes notation and tensor basics that will be used throughout the paper. Section 3 is devoted to the optimization approach for traffic flow data recovery and anomaly diagnosis. Numerical experiments are shown in Section 4, and concluding remarks are drawn in Section 5.

\section{PRELIMINARIES}

\subsection{Notations}

Let $\mathcal{X} \in \mathbb{R}^{I_1 \times \cdots \times I_d}$ denote a $d$th-order tensor, where $x_{i_1 \cdots i_d}$ corresponds to its element. For example, a 3rd-order tensor is represented as $\mathcal{X} \in \mathbb{R}^{I_1 \times I_2 \times I_3}$. Boldface capital letters denote matrices. For instance, an $m \times n$ matrix is written as $\boldsymbol{X} \in \mathbb{R}^{m \times n}$. Boldface lowercase letters represent vectors, such as an $n$-dimensional vector $\boldsymbol{x}_i \in \mathbb{R}^n$, while lowercase letters denote scalars, such as $x_{ij}$. 

For a matrix $\boldsymbol{X}$, its Frobenius norm is defined as $\|\boldsymbol{X}\|_F = \sqrt{\sum_{i,j} x_{ij}^2}$, its $\ell_0$-norm refers to the number of nonzero entries, and its $\ell_{2,0}$-norm represents the number of nonzero rows.  The Frobenius norm of a tensor  $\mathcal{X} \in \mathbb{R}^{I_1 \times \cdots \times I_d}$ is defined as $\|\mathcal{X}\|_F = \sqrt{\sum_{i_1, i_2, \ldots, i_d} x_{i_1 i_2 \cdots i_d}^2}$. The mode-$k$ unfolding of tensor $\mathcal{X}$ is denoted as $\mathcal{X}_{[k]} \in \mathbb{R}^{I_k \times \left( \prod_{l \neq k} I_l \right)}$ for $k = 1, \ldots, d$. Correspondingly, the $\text{fold}_k(\cdot)$ operator folds a matrix into a higher-order tensor along the $k$-th mode, satisfying $\text{fold}_k\left(\mathcal{X}_{[k]}\right) = \mathcal{X}$.

\subsection{Tensor Basics}

\textbf{Definition 1 (Kronecker product \cite{kronecker}).\ } Given matrices $\boldsymbol{A} \in \mathbb{R}^{m \times n}$ and $\boldsymbol{B} \in \mathbb{R}^{p \times q}$, the Kronecker product is defined as:  
\begin{equation*}
\boldsymbol{A} \otimes \boldsymbol{B} := 
\begin{bmatrix}
a_{11}\boldsymbol{B} & a_{12}\boldsymbol{B} & \cdots & a_{1n}\boldsymbol{B} \\
a_{21}\boldsymbol{B} & a_{22}\boldsymbol{B} & \cdots & a_{2n}\boldsymbol{B} \\
\vdots & \vdots & \ddots & \vdots \\
a_{m1}\boldsymbol{B} & a_{m2}\boldsymbol{B} & \cdots & a_{mn}\boldsymbol{B}
\end{bmatrix} \in \mathbb{R}^{(mp) \times (nq)},
\end{equation*}
where the Kronecker product is non-commutative in general.

\textbf{Definition 2 (Mode-$n$ product \cite{Kolda2009}).\ }
Given an \( N \)-th order tensor \( \mathcal{X} \in \mathbb{R}^{I_1 \times I_2 \times \cdots \times I_N} \) and a matrix \( \boldsymbol{A} \in \mathbb{R}^{m \times I_n} \), the mode-\( n \) product is defined as:
\begin{equation*}
    \mathcal{Y} = \mathcal{X} \times_n \boldsymbol{A},
\end{equation*}
where the resulting tensor \( \mathcal{Y} \in \mathbb{R}^{I_1 \times \cdots \times I_{n-1} \times m \times I_{n+1} \times \cdots \times I_N} \) has entries:
\begin{equation*}
    \mathcal{Y}_{i_1, \ldots, i_{n-1}, j, i_{n+1}, \ldots, i_N} = \sum_{i_n=1}^{I_n} \mathcal{X}_{{i_1, \cdots, i_{n-1},i_n, i_{n+1},\cdots, i_N}} \cdot \boldsymbol{A}_{j, i_n}.
\end{equation*}

$\textbf{Definition 3 (Tucker decomposition \cite{Kolda2009}).\ }$
Given $\mathcal{X} \in \mathbb{R}^{I_1 \times \cdots \times I_N}$, the Tucker decomposition of $\mathcal{X}$ can be defined as
\begin{equation*}\begin{array}
{rcl}{\mathcal{X}} & \cong & \sum_{r_{1}=1}^{R_{1}}\cdots\sum_{r_{N}=1}^{R_{N}}g_{r_{1}r_{2}\cdots r_{N}}\left(\boldsymbol{b}_{r_{1}}^{(1)}\circ\boldsymbol{b}_{r_{2}}^{(2)}\circ\cdots\circ\boldsymbol{b}_{r_{N}}^{(N)}\right) \\
 & = & {\mathcal{G}}\times_{1}\boldsymbol{B}^{(1)}\times_{2}\boldsymbol{B}^{(2)}\cdots\times_{N}\boldsymbol{B}^{(N)} \\
 & := & [\![{\mathcal{G}};\{\boldsymbol{B}^{(i)}\}_{i = 
 1}^N]\!],
\end{array}\end{equation*}
where $\mathcal{G} \in \mathbb{R}^{R_1 \times \cdots \times R_N}$ is the core tensor and $\boldsymbol{B}^{(n)} = [\boldsymbol{b}_1^{(n)},\cdots,\boldsymbol{b}_{R_n}^{(n)}]\in \mathbb{R}^{I_n \times R_n}$ are the mode-$n$ factor matrices for $n = 1,\dots,N$.

\section{OPTIMIZATION APPROACH}
\subsection{TSLTO}
For the problem of traffic data completion, we represent the traffic data as a third-order tensor, with road segments as the first dimension, time intervals during the day as the second dimension, and days as the third dimension. The traffic data is then decomposed into low-rank and sparse tensors. Let $\mathcal{Y} \in \mathbb{R}^{D_1 \times D_2 \times D_3}$ denote the observed data, where the set of observed indices is denoted as $\Omega$. The following statistical model is considered:

\begin{equation*}
\mathcal{Y} = \mathcal{W} + \mathcal{R} + \mathcal{\epsilon},
\end{equation*}
where $\mathcal{W} \in \mathbb{R}^{D_1 \times D_2 \times D_3}$ represents the regular traffic data, $\mathcal{R} \in \mathbb{R}^{D_1 \times D_2 \times D_3}$ represents the anomalous traffic data and $\mathcal{\epsilon} \in \mathbb{R}^{D_1 \times D_2 \times D_3}$ represents the error tensor.

Given that regular traffic data exhibits spatio-temporal similarities, including spatial homogeneity and temporal periodicity, $\mathcal{W}$ is low-rank. Thus, we aim to exploit its global low-rank structure through Tucker decomposition:
\begin{equation*}
    \mathcal{W}=[\![\mathcal{G};\{\boldsymbol{U_{i}}\}_{i=1}^{3}]\!].
\end{equation*}
As for local consistency on $\mathcal{W}$, we apply Toeplitz constraints to the factor matrices. Specifically, we consider $\ell_{2,0}$-norm regularization on the following terms:
$\boldsymbol{T}_1\boldsymbol{U}_1,\:\boldsymbol{T}_2\boldsymbol{U}_2\text{ and }\boldsymbol{T}_3\boldsymbol{U}_3$, where $\boldsymbol{T}_1 \in \mathbb{R}^{(D_1-1)\times D_1},\boldsymbol{T}_2\in \mathbb{R}^{(D_2-1)\times D_2},\boldsymbol{T}_3\in \mathbb{R}^{(D_3-1)\times D_3}$ are   Toeplitz matrices:
\begin{equation*}\boldsymbol{T}_i=
\begin{bmatrix}
1 & -1 & 0 & \ldots & 0 & 0 \\
0 & 1 & -1 & \ldots & 0 & 0 \\
 & \vdots&\vdots & \ddots&\vdots &\vdots \\
0 & 0 & 0 & \ldots & 1 & -1
\end{bmatrix} \in\mathbb{R}^{(D_i-1)\times D_i},  \quad i = 1,2,3.\end{equation*}

Besides, as anomalous events occur with low frequency, we impose $\ell_{0}$-norm regularization on $\mathcal{R}$ to enforce sparsity in its structure. Due to the stochastic nature of anomaly events and the fact that such events can affect data within a certain spatiotemporal range, $\mathcal{R}$ exhibits block sparsity. For this, we apply Toeplitz matrix constraints to $\mathcal{R}_{[1]}\in \mathbb{R}^{D_1\times D_2D_3}$, the matricization of $\mathcal{R}$ along its first dimension. Specifically, we impose $\ell_{0}$-norm regularization on the term:
$\boldsymbol{T}_{l}\mathcal{R}_{[1]}\boldsymbol{T}_{r}^{\top}$, where $\boldsymbol{T}_{l}\in \mathbb{R}^{(D_1-1)\times D_1}$, $\boldsymbol{T}_{r}\in \mathbb{R}^{(D_2D_3-1)\times D_2D_3}$.

Let $\Omega$ denotes observed indices, and define $(\cdot)_{\Omega}$ as
\begin{equation*}\left.{(\mathcal{X})_{\Omega}}_{i,j,k} = \left\{\begin{array}{ll}{\mathcal{X}}_{i,j,k},&\text{if}\ (i,j,k) \in \Omega,\\0,&\text{otherwise}.\end{array}\right.\right.\end{equation*}
We  present the optimization problem based on Tucker decomposition:
\begin{equation}
    \begin{aligned}
        \min_{\begin{array}{c}\mathcal{G},\{\boldsymbol{U}_{i}\}_{i=1}^{3},\mathcal{R}\end{array}}&\frac{\beta}{2}\left\|\left([\![\mathcal{G};\{\boldsymbol{U_{i}}\}_{i=1}^{3}]\!]+\mathcal{R}-\mathcal{Y}\right)_{\Omega}\right\|_{F}^{2} + \sum_{i=1}^{3}\lambda_{i}\left\|\boldsymbol{T}_{i}\boldsymbol{U}_{i}\right\|_{2,0} + \mu_{1}\left\|\mathcal{R}\right\|_{0} + \mu_{2}|\left\|\boldsymbol{T}_{l}\mathcal{R}_{[1]}\boldsymbol{T}_{r}^{\top}\right\|_{0} \\
        \text{$s.t.$} \quad 
         &\boldsymbol{U}_{i}^{\top}\boldsymbol{U}_{i} = \boldsymbol{I}_{i}, \quad i=1,2,3,
    \end{aligned}
    \label{orig_model}
\end{equation}
where $\beta>0$, $\mathcal{G}\in\mathbb{R}^{r_1\times r_2\times r_3}$, $\boldsymbol{U}_1\in\mathbb{R}^{D_1\times r_1}$, $\boldsymbol{U}_2\in\mathbb{R}^{D_2\times r_2}$, $\boldsymbol{U}_3\in\mathbb{R}^{D_3\times  r_3}$, $\boldsymbol{R}\in\mathbb{R}^{D_1\times D_2\times D_3}$, $\boldsymbol{T}_l\in\mathbb{R}^{(D_1-1)\times D_1}$ and $\boldsymbol{T}_r\in\mathbb{R}^{(D_2D_3-1)\times D_2D_3}$.  $\boldsymbol{I}_{i}\in \mathbb{R}^{r_i\times r_i}\text{ denotes identity matrix, }  i=1,2,3$. $\mathcal{Y} \in\mathbb{R}^{D_1\times D_2\times D_3}$ is the observed traffic data with $\Omega$.
We know that the recovered data $\mathcal{X}$ is
\begin{equation*}
    \begin{aligned}
       \mathcal{X}=\llbracket \mathcal{G};\{\boldsymbol{U_{i}}\}_{i=1}^{3}\rrbracket+\mathcal{R}.
    \end{aligned}
\end{equation*}
In the proposed optimization model \eqref{orig_model}, the first term is applied to characterize the approximation of the original data at the observed locations. The second term, the $\ell_{2,0}$-norm regulation of $\boldsymbol{T}_{i}\boldsymbol{U}_{i}$, describes the low-rankness of the tensor, which adds constraints to the row sparsity of $\boldsymbol{T}_{i}\boldsymbol{U}_{i}$, specifically the number of distinct rows of factor matrices. The third term is the $\ell_{0}$-norm regulation of the anomaly tensor, which represents the sporadic and stochastic nature of anomalous events, specifically the sparsity of the anomaly tensor. The fourth term, the $\ell_{0}$-norm regulation of $\boldsymbol{T}_{l}\mathcal{R}_{[1]}\boldsymbol{T}_{r}^{\top}$, explains the spatio-temperal consistency of anomalous events by left and right Toeplitz matrix multiplication, revealing the similarity of neighboring abnormal data.

We introduce auxiliary variables $\mathcal{L}$, $\mathcal{X} \in\mathbb{R}^{D_1\times D_2\times D_3}$ and turn the previous model \eqref{orig_model} into
\begin{equation}\begin{aligned}
        \min_{\begin{array}{c}\mathcal{G},\{\boldsymbol{U}_{i}\}_{i=1}^{3},\mathcal{L},\mathcal{R},\mathcal{X}\end{array}}&\frac{\beta}{2}\|\left\llbracket \mathcal{G};\{\boldsymbol{U_{i}}\}_{i=1}^{3}\rrbracket-\mathcal{L}\right\|_{F}^{2} + \sum_{i=1}^{3}\lambda_{i}\left\|\boldsymbol{T}_{i}\boldsymbol{U}_{i}\right\|_{2,0} + \mu_{1}\left\|\mathcal{R}\right\|_{0} + \mu_{2}\left\|\boldsymbol{T}_{l}\mathcal{R}_{[1]}\boldsymbol{T}_{r}^{\top}\right\|_{0} \\
        \text{$s.t.$} \quad 
         &\boldsymbol{U}_{i}^{\top}\boldsymbol{U}_{i} = \boldsymbol{I}_{i}, \quad i=1,2,3, \\
&\mathcal{X} = \mathcal{L}+\mathcal{R},\\
&(\mathcal{X})_{\Omega} = (\mathcal{Y})_{\Omega}.\end{aligned}\end{equation}
The resolution of this problem presents significant challenges, primarily stemming from: (i) the non-convex and discontinuous nature of the $\ell_{2,0}$-norm and the $\ell_0$-norm, (ii) the absence of closed-form solutions for proximal operators associated with linear transformation, (iii) the necessity to satisfy column orthogonality constraints on $U_i$ during the optimization process.

\subsection{ADMM}

Firstly, to solve the non-smooth term and orthogonal constraint efficiently, auxiliary variables $\boldsymbol{Y}_{i} = \boldsymbol{T}_{i}\boldsymbol{U}_{i}$ and $\boldsymbol{Z} = \boldsymbol{T}_{l}\mathcal{R}_{[1]}\boldsymbol{T}_{r}^{\top},\  i = 1,2,3$ are introduced to decouple the linear transformations from non-smooth terms including the $\ell_{2,0}$-norm and the $\ell_{0}$-norm, thereby facilitating the derivation of their proximal operators. The optimization model are reformulated as follows:
\begin{equation}\label{problem12}\begin{aligned}
        \min_{\begin{array}{c}\mathcal{G},\{\boldsymbol{U}_{i}\}_{i=1}^{3},\mathcal{L},\mathcal{R},\mathcal{X},\{\boldsymbol{Y}_{i}\}_{i=1}^{3},\boldsymbol{Z}\end{array}}&\frac{\beta}{2}\|\left\llbracket \mathcal{G};\{\boldsymbol{U_{i}}\}_{i=1}^{3}\rrbracket-\mathcal{L}\right\|_{F}^{2}+\sum_{i=1}^{3}\lambda_{i}\left\|\boldsymbol{Y}_{i}\right\|_{2,0} + \mu_{1}\left\|\mathcal{R}\right\|_{0} + \mu_{2}\left\|\boldsymbol{Z}\right\|_{0} \\
        \text{$s.t.$} \quad & \boldsymbol{Y}_{i} = \boldsymbol{T}_{i}\boldsymbol{U}_{i}, \quad i=1,2,3 \\
        & \boldsymbol{Z} = \boldsymbol{T}_{l}\mathcal{R}_{[1]}\boldsymbol{T}_{r}^{\top} \\
        & \boldsymbol{U}_{i}^{\top}\boldsymbol{U}_{i} = \boldsymbol{I}_{i}, \quad i=1,2,3\\
&\mathcal{X} = \mathcal{L}+\mathcal{R}\\
&(\mathcal{X})_{\Omega} = (\mathcal{Y})_{\Omega},
    \end{aligned}\end{equation}
Taking the advantage of the block-separable structure from the above resulting optimization model, we adopt the Alternating Direction Method of Multipliers (ADMM) for numerical resolution. As for the constraints $\boldsymbol{U}_{i}^{\top}\boldsymbol{U}_{i} = \boldsymbol{I}_{i} \quad i = 1,2,3$, we encode the column orthogonality requirement via an indicator function within the augmented Lagrangian framework instead of handling this constraint through conventional methods. This formulation guarantees that $\boldsymbol{U_i}$ satisfies the orthogonality of columns throughout the optimization process. Let $\ D_i =dim(\mathcal{Y})(i)$, $ St(D_i,r_i)$ is defined as $\{\boldsymbol{U}_i\in\mathbb{R}^{D_i\times r_i}:\boldsymbol{U}_i^{\top}\boldsymbol{U}_i=I_{i}\}$, the indicator function can be expressed as:
\begin{equation*}\left.\delta_{St(D_i,r_i)}(\boldsymbol{U}_{i})=\left\{\begin{array}{ll}\mathbf{0},&\text{if }\ \boldsymbol{U}_{i} \in St(D_i,r_i),\\+\infty,&\text{otherwise}.\end{array}\right.\right.\end{equation*}All remaining equality constraints have been associated with dedicated Lagrange multipliers, with corresponding linear and quadratic terms incorporated into the augmented Lagrangian function. The complete augmented Lagrangian for the optimization problem is formulated as
\begin{equation*}\begin{aligned}&\mathscr{L}(\mathcal{X},\mathcal{G},\{\boldsymbol{U}_{i}\}_{i=1}^{3},\mathcal{R},\mathcal{L},\{\boldsymbol{Y}_{i}\}_{i=1}^{3},\boldsymbol{Z},\{\boldsymbol{V}_{i}\}_{i=1}^{3},\boldsymbol{W},\mathcal{P})\\=&
    \frac{\beta}{2}\|\left \llbracket \mathcal{G};\{\boldsymbol{U_{i}}\}_{i=1}^{3}\rrbracket-\mathcal{L}\right\|_{F}^{2}+\sum_{i=1}^{3}\lambda_{i}\|\boldsymbol{Y}_{i}\|_{2,0}+\mu_{1}\|\mathcal{R}\|_{0}
    +\mu_{2}\|\boldsymbol{Z}\|_{0}\\+&\langle \boldsymbol{Z}-\boldsymbol{T}_{l}\mathcal{R}_{[1]}\boldsymbol{T}_{r}^{\top},\boldsymbol{W}\rangle
    +\frac{\gamma}{2}\|\boldsymbol{Z}-\boldsymbol{T}_{l}\mathcal{R}_{[1]}\boldsymbol{T}_{r}^{\top}\|_F^{2}+\sum_{i=1}^{3}(\langle \boldsymbol{Y}_{i}-\boldsymbol{T}_{i}\boldsymbol{U}_{i} ,\boldsymbol{V}_{i}\rangle\\+&\frac{\alpha_{i}}{2}\|\boldsymbol{Y}_{i}-\boldsymbol{T}_{i}\boldsymbol{U}_{i}\|_{F}^{2})+\langle \mathcal{X} - \mathcal{L} - \mathcal{R}, \mathcal{P}\rangle + \frac{s}{2}\|\mathcal{X} - \mathcal{L} - \mathcal{R}\|_F^{2}+\sum_{i=1}^{3}\delta_{St(D_i,r_i)}(\boldsymbol{U}_{i}).
\end{aligned}\end{equation*}

Next, we update each variable iteratively:
\begin{equation*}\begin{cases}
\mathcal{X}^{k+1}&:= \arg\min_{\mathcal{X}}\mathscr{L}(\mathcal{X},\mathcal{G}^k,\{\boldsymbol{U}_{i}^k\}_{i=1}^{3},\mathcal{R}^k,\mathcal{L}^k,\{\boldsymbol{Y}_{i}^k\}_{i=1}^{3},\boldsymbol{Z}^k,\{\boldsymbol{V}_{i}^k\}_{i=1}^{3},\boldsymbol{W}^k,\mathcal{P}^k),\\
\mathcal{G}^{k+1}&:= \arg\min_{\mathcal{G}}\mathscr{L}(\mathcal{X}^{k+1},\mathcal{G},\{\boldsymbol{U}_{i}^k\}_{i=1}^{3},\mathcal{R}^k,\mathcal{L}^k,\{\boldsymbol{Y}_{i}^k\}_{i=1}^{3},\boldsymbol{Z}^k,\{\boldsymbol{V}_{i}^k\}_{i=1}^{3},\boldsymbol{W}^k,\mathcal{P}^k),\\
{\boldsymbol{U}_1}^{k+1}&:= \arg\min_{\boldsymbol{U}_1}\mathscr{L}(\mathcal{X}^{k+1},\mathcal{G}^{k+1},\boldsymbol{U}_{1},\boldsymbol{U}_{2}^{k},\boldsymbol{U}_{3}^{k},\mathcal{R}^k,\mathcal{L}^k,\{\boldsymbol{Y}_{i}^k\}_{i=1}^{3},\boldsymbol{Z}^k,\{\boldsymbol{V}_{i}^k\}_{i=1}^{3},\boldsymbol{W}^k,\mathcal{P}^k),\\
{\boldsymbol{U}_2}^{k+1}&:= \arg\min_{\boldsymbol{U}_2}\mathscr{L}(\mathcal{X}^{k+1},\mathcal{G}^{k+1},\boldsymbol{U}_{1}^{k+1},\boldsymbol{U}_{2}^{},\boldsymbol{U}_{3}^{k},\mathcal{R}^k,\mathcal{L}^k,\{\boldsymbol{Y}_{i}^k\}_{i=1}^{3},\boldsymbol{Z}^k,\{\boldsymbol{V}_{i}^k\}_{i=1}^{3},\boldsymbol{W}^k,\mathcal{P}^k),\\
{\boldsymbol{U}_3}^{k+1}&:= \arg\min_{\boldsymbol{U}_3}\mathscr{L}(\mathcal{X}^{k+1},\mathcal{G}^{k+1},\boldsymbol{U}_{1}^{k+1},\boldsymbol{U}_{2}^{k+1},\boldsymbol{U}_{3}^{},\mathcal{R}^k,\mathcal{L}^k,\{\boldsymbol{Y}_{i}^k\}_{i=1}^{3},\boldsymbol{Z}^k,\{\boldsymbol{V}_{i}^k\}_{i=1}^{3},\boldsymbol{W}^k,\mathcal{P}^k),\\
{R}^{k+1} &:= arg \min_{\mathcal{R}} \mathscr{L}(\mathcal{X}^{k+1},\mathcal{G}^{k+1},\{\boldsymbol{U}_{i}^{k+1}\}_{i=1}^{3},\mathcal{R},\mathcal{L}^k,\{\boldsymbol{Y}_{i}^k\}_{i=1}^{3},\boldsymbol{Z}^k,\{\boldsymbol{V}_{i}^k\}_{i=1}^{3},\boldsymbol{W}^k,\mathcal{P}^k),\\
\mathcal{L}^{k+1} &:=arg \min_{\mathcal{L}} \mathscr{L}(\mathcal{X}^{k+1},\mathcal{G}^{k+1},\{\boldsymbol{U}_{i}^{k+1}\}_{i=1}^{3},\mathcal{R}^{k+1},\mathcal{L},\{\boldsymbol{Y}_{i}^k\}_{i=1}^{3},\boldsymbol{Z}^k,\{\boldsymbol{V}_{i}^k\}_{i=1}^{3},\boldsymbol{W}^k,\mathcal{P}^k),\\
\boldsymbol{Y}_i^{k+1} &:=arg \min_{\boldsymbol{Y}_i} \mathscr{L}(\mathcal{X}^{k+1},\mathcal{G}^{k+1},\{\boldsymbol{U}_{i}^{k+1}\}_{i=1}^{3},\mathcal{R}^{k+1},\mathcal{L}^{k+1},\boldsymbol{Y}_{i},\boldsymbol{Z}^k,\{\boldsymbol{V}_{i}^k\}_{i=1}^{3},\boldsymbol{W}^k,\mathcal{P}^k),\quad i=1,2,3\\
\boldsymbol{Z}^{k+1} &:=arg \min_{\boldsymbol{Z}} \mathscr{L}(\mathcal{X}^{k+1},\mathcal{G}^{k+1},\{\boldsymbol{U}_{i}^{k+1}\}_{i=1}^{3},\mathcal{R}^{k+1},\mathcal{L}^{k+1},\{\boldsymbol{Y}_{i}^{k+1}\}_{i=1}^{3},\boldsymbol{Z},\{\boldsymbol{V}_{i}^k\}_{i=1}^{3},\boldsymbol{W}^k,\mathcal{P}^k),\\
\boldsymbol{V}_i^{k+1}&:=\boldsymbol{V}_i^k+\alpha_i(\boldsymbol{Y}_{i}^{k+1}-\boldsymbol{T}_{i}\boldsymbol{U}^{k+1}_{i}),\quad i = 1,2,3,\\
\boldsymbol{W}^{k+1} &:= \boldsymbol{W}^k+\gamma(\boldsymbol{Z}^{k+1}-\boldsymbol{T}_{l}\mathcal{R}_{[1]}^{k+1}\boldsymbol{T}_{r}^{\top}),\\
\mathcal{P}^{k+1} &:= \mathcal{P}^{k} + s(\mathcal{X}^{k+1} -\mathcal{L}^{k+1}- \mathcal{R}^{k+1}).
\end{cases}\end{equation*}

\textbf{Update $\mathcal{X}$}: A closed-form solution via projections related to the index set $\Omega$ takes the form of 
\begin{equation}
\begin{aligned}
\mathcal{X}^{k+1} &:= \arg\min_{\mathcal{X}} \left\{ \langle \mathcal{X} - \mathcal{L}^k - \mathcal{R}^k, \mathcal{P}^k \rangle + \frac{s}{2} \|\mathcal{X} - \mathcal{L}^k - \mathcal{R}^k\|_F^2 \;\Big|\; (\mathcal{X})_{\Omega} = (\mathcal{Y})_{\Omega} \right\} \\
    &=(\mathcal{Y})_{\Omega} + (\mathcal{L}^k+\mathcal{R}^k - \frac{1}{s}\mathcal{P}^k)_{\bar{\Omega}}.
\end{aligned}
\label{eq:update_X}
\end{equation}

\textbf{Update $\mathcal{G}$}: The closed-form solution resulting from the first-order optimality gives us
%
\begin{equation}
\begin{aligned}
    \mathcal{G}^{k+1} &:= arg \min_{\mathcal{G}}\frac{\beta}{2}\|\left\llbracket\mathcal{G};\{\boldsymbol{U_{i}^k}\}_{i=1}^{3}\rrbracket-\mathcal{L}^k\right\|_{F}^{2}\\
    &=\mathcal{L}^k\times_{1}{\boldsymbol{U}_{1}^{k}}^{\top}\times_{2}{\boldsymbol{U}_{2}^{k}}^{\top}\times_{3}{\boldsymbol{U}_{3}^{k}}^{\top}.
\end{aligned}
\label{eq:update_G}
\end{equation}

\textbf{Update $\boldsymbol{U}$}:
The corresponding subproblems quadratic programming over Stiefel manifolds, see as below, which will be approximately solved by
gradient-type manifold optimization approach with a curvilinear search method incorporating Barzilai-Borwein (BB) step size adaptation \cite{OptM-Wen-Yin-2010}.
\begin{equation}
\begin{aligned}
    {\boldsymbol{U}_1}^{k+1} &:= arg \min_{\boldsymbol{U}_1} \frac{\beta}{2}\|\mathcal{G}^{k+1}\times_{1}\boldsymbol{U}_{1}\times_{2}\boldsymbol{U}_{2}^{k}\times_{3}\boldsymbol{U}_{3}^k-\mathcal{L}^k\|_{F}^{2} \\&+ \langle \boldsymbol{Y}_{1}^k-\boldsymbol{T}_{1}\boldsymbol{U}_{1} ,\boldsymbol{V}_{1}^k\rangle +\frac{\alpha_{1}}{2}||\boldsymbol{Y}_{1}^k-\boldsymbol{T}_{1}\boldsymbol{U}_{1}||_{F}^{2} + \delta_{St(D_1,r_1)}(\boldsymbol{U}_{1}),
\end{aligned}
\label{eq:update_U1}
\end{equation}
\begin{equation}
\begin{aligned}
    {\boldsymbol{U}_2}^{k+1} &:= arg \min_{\boldsymbol{U}_2} \frac{\beta}{2}\|\mathcal{G}^{k+1}\times_{1}\boldsymbol{U}_{1}^{k+1}\times_{2}\boldsymbol{U}_{2}\times_{3}\boldsymbol{U}_{3}^k-\mathcal{L}^k\|_{F}^{2} \\&+ \langle \boldsymbol{Y}_{2}^k-\boldsymbol{T}_{2}\boldsymbol{U}_{2} ,\boldsymbol{V}_{2}^k\rangle +\frac{\alpha_{2}}{2}||\boldsymbol{Y}_{2}^k-\boldsymbol{T}_{2}\boldsymbol{U}_{2}||_{F}^{2} + \delta_{St(D_2,r_2)}(\boldsymbol{U}_{2}),
\end{aligned}
\label{eq:update_U2}
\end{equation}
\begin{equation}
\begin{aligned}
    {\boldsymbol{U}_3}^{k+1} &:= arg \min_{\boldsymbol{U}_3} \frac{\beta}{2}\|\mathcal{G}^{k+1}\times_{1}\boldsymbol{U}_{1}^{k+1}\times_{2}\boldsymbol{U}_{2}^{k+1}\times_{3}\boldsymbol{U}_{3}-\mathcal{L}^k\|_{F}^{2} \\&+ \langle \boldsymbol{Y}_{3}^k-\boldsymbol{T}_{3}\boldsymbol{U}_{3} ,\boldsymbol{V}_{3}^k\rangle +\frac{\alpha_{3}}{2}||\boldsymbol{Y}_{3}^k-\boldsymbol{T}_{3}\boldsymbol{U}_{3}||_{F}^{2} + \delta_{St(D_3,r_3)}(\boldsymbol{U}_{3}).
\end{aligned}
\label{eq:update_U3}
\end{equation}
The gradients with respect to $\boldsymbol{U}_{i}$  of 
$
    f_\beta([\![\mathcal{G};\boldsymbol{U}_1,\boldsymbol{U}_2,\boldsymbol{U}_3]\!],\mathcal{L}) = \frac{\beta}{2}\left\|[\![\mathcal{G};\{\boldsymbol{U}_{i}\}_{i=1}^{3}]\!]-\mathcal{L}\right\|_{F}^{2}
$ are calculated by the following equations:
\begin{equation*}
   \frac{\partial f_\beta}{\partial \boldsymbol{U}_{1}} =
   \beta ([\![\mathcal{G};\{\boldsymbol{U}_{i}\}_{i=1}^{3}]\!]-\mathcal{L})_{(1)}(\boldsymbol{U}_{3} \otimes \boldsymbol{U}_{2})\mathcal{G}_{(1)}^\top,
\end{equation*}
\begin{equation*}
    \frac{\partial f_\beta}{\partial \boldsymbol{U}_{2}} =\beta ([\![\mathcal{G};\{\boldsymbol{U}_{i}\}_{i=1}^{3}]\!]-\mathcal{L})_{(2)}(\boldsymbol{U}_{3} \otimes \boldsymbol{U}_{1})\mathcal{G}_{(2)}^\top,
\end{equation*}
\begin{equation*}
    \frac{\partial f_\beta}{\partial \boldsymbol{U}_{3}} =\beta ([\![\mathcal{G};\{\boldsymbol{U}_{i}\}_{i=1}^{3}]\!]-\mathcal{L})_{(3)}(\boldsymbol{U}_{2} \otimes \boldsymbol{U}_{1})\mathcal{G}_{(3)}^\top.
\end{equation*}

\textbf{Update $\mathcal{R}$}:
To deal with the absence of a closed-form proximal operator for the $\ell_{0}$-norm under linear transformations, we employ proximal gradient method to update $R$ approximately.
\begin{equation}
    \begin{aligned}\mathcal{R}^{k+1} &= arg \min_{\mathcal{R}} \mathscr{L}(\mathcal{G}^{k+1},\{\boldsymbol{U}_{i}^{k+1}\}_{i=1}^{3},\mathcal{R},\{\boldsymbol{Y}_{i}^k\}_{i=1}^{3},\boldsymbol{Z}^k,\{\boldsymbol{V}^k_{i}\}_{i=1}^{3},\boldsymbol{W}^k))\\
     &=arg \min_{\mathcal{R}}\mu_{1}||\mathcal{R}||_{0}+\langle \boldsymbol{Z}^{k}-\boldsymbol{T}_{l}\mathcal{R}_{[1]}\boldsymbol{T}_{r}^{\top},\boldsymbol{W}^{k}\rangle+\frac{\gamma}{2}\|\boldsymbol{Z}^{k}-\boldsymbol{T}_{l}\mathcal{R}_{[1]}\boldsymbol{T}_{r}^{\top}\|_{F}^{2} \\&+\langle \mathcal{X}^k-\mathcal{L}^k-\mathcal{R}, \mathcal{P}^k \rangle + \frac{s}{2}\|\mathcal{X}^k-\mathcal{L}^k-\mathcal{R}\|_{F}^{2}
\\& \approx\mathbf{prox}_{\lambda^k \mathcal{\mu}_1\|\cdot\|_{0}} \Bigg(  \mathcal{R}^k - \lambda^k \Big[ 
 - \text{fold}_k(\boldsymbol{T}_{l}^{\top} \boldsymbol{W}^k \boldsymbol{T}_{r})
  \\&- \gamma \ \text{fold}_k(\boldsymbol{T}_{l}^{\top}(\boldsymbol{Z}^k - \boldsymbol{T}_{l} \mathcal{R}_{[1]}^k \boldsymbol{T}_{r}^{\top}) \boldsymbol{T}_{r})-\mathcal{P}^k-s(\mathcal{X}^k-\mathcal{L}^k-\mathcal{R}^k) \Big] \Bigg),
\end{aligned}
\label{eq:update_R}
\end{equation}
where $\lambda^k$ is obtained by line search. 

\textbf{Update $\mathcal{L}$}: First-order optimality yields
\begin{equation}
    \begin{aligned}
        \mathcal{L}^{k+1} &= arg\min_{\mathcal{L}} \frac{\beta}{2}\|\left\llbracket\mathcal{G}^{k+1};\{\boldsymbol{U}_{i}^{k+1}\}_{i=1}^{3}\rrbracket-\mathcal{L}\right\|_{F}^{2}\\&+\langle \mathcal{X}^{k+1} - \mathcal{L} - \mathcal{R}^{k+1}, \mathcal{P}^k\rangle + \frac{s}{2}\|\mathcal{X}^{k+1} - \mathcal{L} - \mathcal{R}^{k+1}\|_F^{2}\\
        &= \frac{\beta}{\beta + s}[\![\mathcal{G}^{k+1};\{\boldsymbol{U}_{i}^{k+1}\}_{i=1}^{3}]\!]+\frac{s}{\beta + s}(\mathcal{X}^{k+1} - \mathcal{R}^{k+1})+\frac{1}{\beta + s}\mathcal{P}^k.
    \end{aligned}
    \label{eq:update_L}
\end{equation}

\textbf{Update $\boldsymbol{Y}$}: The closed-form solution achieves via the underlying proximal operator:
\begin{equation}\begin{aligned}\boldsymbol{Y}_i^{k+1} &= arg \min_{\boldsymbol{Y}_i} \mathscr{L}(\mathcal{G}^{k+1},\{\boldsymbol{U}_{i}^{k+1}\}_{i=1}^{3},\mathcal{R}^{k+1},\boldsymbol{Y}_{i},\boldsymbol{Z}^k,\{\boldsymbol{V}^k_{i}\}_{i=1}^{3},\boldsymbol{W}^k)\\&= arg \min_{\boldsymbol{Y}_i}\lambda_{i}||\boldsymbol{Y}_{i}||_{2,0}+\langle \boldsymbol{Y}_{i}-\boldsymbol{T}_{i}\boldsymbol{U}^{k+1}_{i},\boldsymbol{V}^{k}_{i}\rangle+\frac{\alpha_{i}}{2}||\boldsymbol{Y}_{i}-\boldsymbol{T}_{i}\boldsymbol{U}^{k+1}_{i}||_{F}^{2}\\&= arg \min_{\boldsymbol{Y}_i}\frac{\lambda_i}{\alpha_{i}}||\boldsymbol{Y}_{i}||_{2,0}+\frac{1}{\alpha_i} \langle \boldsymbol{Y}_{i}-\boldsymbol{T}_{i}\boldsymbol{U}^{k+1}_{i},\boldsymbol{V}^{k}_{i}\rangle+ \frac{1}{2}||\boldsymbol{Y}_{i}-\boldsymbol{T}_{i}\boldsymbol{U}^{k+1}_{i}||_{F}^{2}\\&=\text{prox}_{\frac{\lambda_i}{\alpha_i} \|\cdot\|_{2,0}}(\boldsymbol{T_i} \boldsymbol{U_i}^{k+1}-\frac{1}{\alpha_i} \boldsymbol{V_i}^k), \quad i=1,2,3,\end{aligned}\label{eq:update_Y}\end{equation}
where the group hard-thresholding operator $\text{prox}_{\frac{\lambda_i}{\alpha_i} \|\cdot\|_{2,0}}(\cdot)$ can be chosen row-wisely as:
\begin{equation*}\left.\boldsymbol{y}_i^{j}=\left\{\begin{array}{ll}\mathbf{0},&||\boldsymbol{x}_i^{j}||_{2}^{2}\leq\frac{2\lambda_i}{\alpha_i},\\\boldsymbol{x}_i^{j},&||\boldsymbol{x}_i^{j}||_{2}^{2}>\frac{2\lambda_i}{\alpha_i},\end{array}\right.\right.\end{equation*}
with $\boldsymbol{y_i}^{j}$ the $j$-th row of $\boldsymbol{Y}_i$ and $\boldsymbol{x}_i^{j}$  the $j$-th row of $\boldsymbol{T_i} \boldsymbol{U_i}^{k+1}-\frac{1}{\alpha_i} \boldsymbol{V_i}^k$.

\textbf{Update $\boldsymbol{Z}$}:
\begin{equation}
\begin{aligned}
    \boldsymbol{Z}^{k+1} &= arg \min_{\boldsymbol{Z}} \mathscr{L}(\mathcal{G}^{k+1},\{\boldsymbol{U}_{i}^{k+1}\}_{i=1}^{3},\mathcal{R}^{k+1},\boldsymbol{Y}_{i}^{k+1},\boldsymbol{Z},\{\boldsymbol{V}^k_{i}\}_{i=1}^{3},\boldsymbol{W}^k)\\&= arg \min_{\boldsymbol{Z}}\mu_{2}||\boldsymbol{Z}||_{0}+\langle \boldsymbol{Z}-\boldsymbol{T}_{l} \mathcal{R}_{[1]}^{k+1}\boldsymbol{T}_{r}^{\top},\boldsymbol{W}^k\rangle+\frac{\gamma}{2}\|\boldsymbol{Z}-\boldsymbol{T}_{l} \mathcal{R}_{[1]}^{k+1}\boldsymbol{T}_{r}^{\top}\|_F^{2}\\&= arg \min_{\boldsymbol{Z}}\frac{\mu_2}{\gamma}||\boldsymbol{Z}||_{0}+\frac{1}{\gamma} \langle \boldsymbol{Z}-\boldsymbol{T}_{l}\mathcal{R}_{[1]}^{k+1}\boldsymbol{T}_{r}^{\top},\boldsymbol{W}^k\rangle+\frac{1}{2}\|\boldsymbol{Z}-\boldsymbol{T}_{l}\mathcal{R}_{[1]}^{k+1}\boldsymbol{T}_{r}^{\top}\|_F^{2}\\&=\mathbf{prox}_{\frac{\mu_2}{\gamma}\|\cdot\|_{0}}(\boldsymbol{T}_{l}\mathcal{R}_{[1]}^{k+1}\boldsymbol{T}_{r}^{\top}-\frac{1}{\gamma}\boldsymbol{W}^{k}).
\end{aligned}
\label{eq:update_Z}
\end{equation}
Similarly, this subproblem can be directly reduced to computing the proximal operator of the $\ell_0$-norm:
\begin{equation*}\left.{z}_{i,j}=\left\{\begin{array}{ll}\mathbf{0},&{x}_{i,j}^2 \leq\frac{2\mu_2}{\gamma},\\{x}_{i,j},&{x}_{i,j}^2>\frac{2\mu_2}{\gamma},\end{array}\right.\right.\end{equation*}
where ${z}_{i,j}$ denotes the element in the i-th row and j-th column of $\boldsymbol{Z}$, ${x}_{i,j}$ denotes the element in the i-th row and j-th column of $\boldsymbol{T}_{l}\mathcal{R}_{[1]}^{k+1}\boldsymbol{T}_{r}^{\top}-\frac{1}{\gamma}\boldsymbol{W}^{k}$.

\textbf{Update Lagrange Multipliers $\{\boldsymbol{V}_{i}\}_{i=1}^{3}$,$\boldsymbol{W}$,$\mathcal{P}$}:
\begin{equation}
    \boldsymbol{V}_i^{k+1}=\boldsymbol{V}_i^k+\alpha_i(\boldsymbol{Y}_{i}^{k+1}-\boldsymbol{T}_{i}\boldsymbol{U}^{k+1}_{i}), \quad i = 1,2,3,
\label{eq:update_V}
\end{equation}
\begin{equation}
    \boldsymbol{W}^{k+1} = \boldsymbol{W}^k+\gamma(\boldsymbol{Z}^{k+1}-\boldsymbol{T}_{l}\mathcal{R}_{[1]}^{k+1}\boldsymbol{T}_{r}^{\top}),
\label{eq:update_W}
\end{equation}
\begin{equation}
    \mathcal{P}^{k+1} = \mathcal{P}^{k} + s(\mathcal{X}^{k+1} -\mathcal{L}^{k+1}- \mathcal{R}^{k+1}).
    \label{eq:update_P}
\end{equation}

The algorithmic framework of the ADMM for solving 
the proposed sparse low-rank high-order tensor optimization model \eqref{problem12} is listed in Algorithm 1.

\begin{algorithm}[ht]
\setcounter{AlgoLine}{0}
\caption{ADMM for problem \eqref{problem12} }
\KwIn{The collected traffic data tensor $\mathcal{Y}$, the index set of the observed data $\Omega$, tucker size $[r_1,r_2,r_3]$,\\$\alpha$, $\beta$, $\gamma$, $\lambda$, $\mu$, $s$, $\delta$, $\eta$, $c$, $\tau$, $\lambda_0$, $\rho$, $\varepsilon$.}
\KwOut{The recovered tensor $\mathcal{X}$, Low-Rank tensor $\mathcal{L}$ and the block-sparse anomaly tensor $\mathcal{R}$.}
\SetKwInOut{Input}{}
\SetKwInOut{Output}{}

\For{k = 1 \textbf{to} $k_{max}$}{
    \textbf{Update} $\mathcal{X}$ \textbf{by} \eqref{eq:update_X}. \;
    
    \textbf{Update} $\mathcal{G}$ \textbf{by} \eqref{eq:update_G}. \;
    
    \textbf{for} i = 1 \textbf{to} 3,
        \textbf{Update} $\boldsymbol{U_i}$ \textbf{by} \eqref{eq:update_U1} \eqref{eq:update_U2} \eqref{eq:update_U3}. \;
    
    \textbf{Update} $\mathcal{R}$ \textbf{by} \eqref{eq:update_R}. \;

    \textbf{Update} $\mathcal{L}$ \textbf{by} \eqref{eq:update_L}. \;
    
    \textbf{for} i = 1 \textbf{to} 3,
        \textbf{Update} $\boldsymbol{Y_i}$ \textbf{by} \eqref{eq:update_Y}. \;
    
    \textbf{Update} $\boldsymbol{Z}$ \textbf{by} \eqref{eq:update_Z}. \;
    
    \textbf{for} i = 1 \textbf{to} 3,
        \textbf{Update} $\boldsymbol{V_i}$ \textbf{by} \eqref{eq:update_V}. \;
    
    \textbf{Update} $\boldsymbol{W}$ \textbf{by} \eqref{eq:update_W}. \;

    \textbf{Update} $\mathcal{P}$ \textbf{by} \eqref{eq:update_P}. \;
    
    \If{k \textgreater 1 \textbf{and} $\frac{||\mathcal{X}_{pre} - \mathcal{X}||}{||\mathcal{X}_{pre}||}$ \textbf{and}  $\frac{||\mathcal{G}_{pre} - \mathcal{G}||}{||\mathcal{G}_{pre}||}$ \textbf{and}  $\frac{||\mathcal{L}_{pre} - \mathcal{L}||}{||\mathcal{L}_{pre}||}$ \textbf{and}  $\frac{||\mathcal{R}_{pre} - \mathcal{R}||}{||\mathcal{R}_{pre}||}$ \textless $\epsilon$}{
        \textbf{break.} \;
    }
    $\alpha$ = $\alpha$*1.15, $\gamma$ = $\gamma$*1.15, $s$ = $s$*1.15. \;
}
\end{algorithm}

\section{NUMERICAL EXPERIMENTS}
All experiments were performed on a MacBook Pro (2021) equipped with an Apple M1 Max chip (12-core CPU, 32-core GPU, 3.5 GHz) and 32GB of RAM, running macOS Sequoia 15.2. The numerical simulations were executed in MATLAB R2023b with the default multi-threading configuration. No GPU acceleration was  utilized\footnote{Matlab codes are available at \url{https://github.com/TSLTO2025/TSLTO}}.
\subsection{Evaluation Criteria}

The data imputation performance of our model is evaluated by Root Mean Square Error (RMSE), Mean Absolute Percentage Error (MAPE) and Mean Absolute Error (MAE), which are respectively defined as:
\begin{equation*}
    \text{RMSE} = \sqrt{\frac{1}{n} \sum_{i=1}^{n} (x_i - \hat{x}_i)^2},
\end{equation*}
\begin{equation*}
    \text{MAPE} = \frac{1}{n} \sum_{i=1}^{n} \left| \frac{x_i - \hat{x}_i}{x_i} \right| \times 100\%,
\end{equation*}
\begin{equation*}
    \text{MAE} = \frac{1}{n} \sum_{i=1}^{n} |x_i - \hat{x}_i|,
\end{equation*}
where $x_i$ denotes the ground truth value and $\hat{x_i}$
represents the imputed value. Lower RMSE, MAPE and MAE values indicate superior tensor completion performance.

We determines the validity of anomaly detection by binary classification based on spatiotemporal overlap: a positive detection is registered if the identified anomaly intersects with the ground truth in spatial or temporal dimensions; otherwise, it is marked as negative. To evaluate its accuracy, we introduce Precision, Recall, and the F1 Score:
\begin{equation*}
            \text{Precision} = \frac{TP}{TP + FP},
\end{equation*}
\begin{equation*}
            \text{Recall} = \frac{TP}{TP + FN},
\end{equation*}
\begin{equation*}
            \text{F1\ Score} = \frac{2\text{Precision} \text{Recall}}{\text{Precision} + \text{Recall}},            
\end{equation*}
where True Positives (TP) denote the number of correctly identified anomalous instances, False Positives (FP) represent the count of normal instances erroneously classified as anomalous, and False Negatives (FN) indicate the quantity of undetected true anomalies. These three metrics range between 0 and 1, with higher values indicating superior anomaly diagnosis capability.

\subsection{Dataset Description}
\subsubsection{Synthetic Dataset}
We evaluate our model's dual capability in simultaneous data imputation and anomaly diagnosis using a synthetic dataset of dimensions $50 \times 50 \times 50$. The low-rank component is synthesized based on Tucker decomposition with a core tensor $\mathcal{G} \in \mathbb{R}^{3 \times 3 \times 3}$ containing uniformly distributed values $U(0, 100)$, and three factor matrices $\boldsymbol{U_i} \in \mathbb{R}^{50 \times 3}$ ($i = 1,2,3$) with 2 distinctive rows in each matrix contrasting the remaining 48 conventional rows. As for anomaly synthesis, according to central-limit theorem, we suppose it follows a normal distribution $N(4,0.01)$, distributed as 50 discrete blocks of size $2 \times 125$ along the mode-1 matricization. This configuration yields a controlled anomaly contamination ratio of 10\% within the tensor structure. Figure~\ref{moni} illustrates the 12th frontal slice along the mode-3 orientation:
\begin{figure}[H]
    \centering
    \includegraphics[width=0.75\linewidth]{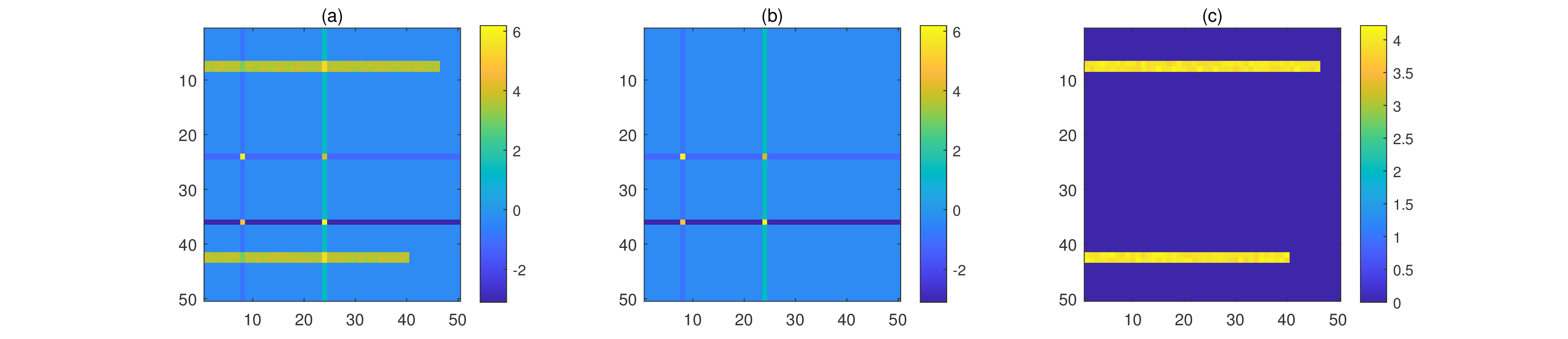}
    \captionsetup{font=footnotesize}
    \caption{Synthetic data slice. Subfigure (a) displays the original synthetic tensor, subfigure (b) visualizes the low-rank component, and subfigure (c) highlights the sparse anomaly component.}
    \label{moni}
\end{figure}

\subsubsection{Real-world Datasets}
We employ a real-world traffic speed dataset for numerical validation: the Guangzhou Urban Traffic Speed Dataset (GUANGZHOU)\footnote{\url{https://zenodo.org/records/1205229}} 
to verify the practicality and efficiency of our proposed tensor model and algorithm in traffic flow data recovery.

The GUANGZHOU dataset comprises urban traffic speed data from 214 road segments recorded at 10-minute intervals over 61 days in 2016. Configured as a third-order tensor with dimensions 214 (road segments) × 144 (daily time intervals) × 61 (days), the dataset contains 1,879,776 entries. Missing values are encoded as zeros, yielding 1,855,589 valid observations with a missing rate of 1.29\%.


\subsection{Parameter Tuning and Evaluation Results}

\begin{table}[ht]
\centering
\caption{Anomaly Diagnosis and Data Completion Evaluation on simulated data}
\label{recall}
\begin{threeparttable}
\begin{tabular}{@{}ccccccccc@{}}
\toprule
Missing rate     & 0.1     & 0.2     & 0.3     & 0.4     & 0.5     & 0.6     & 0.7     & 0.8      \\
\midrule
$\mu_1$          & 5       & 5       & 5       & 5       & 5       & 5       & 5       & 5        \\
$\mu_2$          & 10      & 15      & 20      & 20      & 25      & 35      & 50      & 60       \\
$\beta$          & 300     & 450     & 450     & 400     & 450     & 450     & 450     & 370      \\
$\gamma$         & 10      & 5       & 10      & 10      & 10      & 10      & 10      & 10       \\
$\lambda_i$\tnote{1}      & 1.2     & 1.2     & 1.2     & 1.2     & 1.2       & 1.2     & 1.2     & 1.2      \\
\midrule
Total iterations & 65737   & 38018   & 35001   & 14559   & 35002   & 32062   & 63540   & 47210    \\
Precision        & 0.914   & 0.694   & 0.857   & 0.846   & 0.897   & 0.792   & 0.0.801   & 0.826    \\
Recall           & 0.958   & 0.956   & 0.931   & 0.871   & 0.814   & 0.772   & 0.715   & 0.652    \\
F1 Score         & 0.936   & 0.805   & 0.892   & 0.858   & 0.853   & 0.782   & 0.755   & 0.729    \\
MAPE             & 5.924\% & 3.941\% & 6.506\% & 6.168\% & 6.923\% & 9.328\% & 8.220\% & 11.229\% \\
RMSE             & 0.428   & 0.445   & 0.392   & 0.407   & 0.448   & 0.455   & 0.456   & 0.470    \\
MAE              & 0.035   & 0.039   & 0.037   & 0.037   & 0.039   & 0.049   & 0.047   & 0.060   \\
\bottomrule
\end{tabular}
\begin{tablenotes}
\footnotesize            
    \item[1] If $\lambda_i$ is 1.2, then $\lambda$ is [1.2,1.2,1.2].
    \item[2] Parameters not mentioned above: $\alpha$ = [100,100,100], k = 0.3, max iteration for inner loop is set to 60, and max iteration for outer loop is set to 200,000.
\end{tablenotes}
\end{threeparttable}
\end{table}

We conducted parameter tuning on simulated data, sequentially adjusting the objective function coefficients of the original optimization problem in the order of $\beta$, $\mu_1$, $\mu_2$, and $\gamma$. Our experimental observations revealed that increasing $\mu_2$ enhances the $\mathrm{F1}$-Score performance, and results in deterioration of the evaluation metrics $\mathrm{MAPE}$, $\mathrm{RMSE}$, and $\mathrm{MAE}$. Besides, appropriately increasing $\beta$ while decreasing $\gamma$ was found to improve the evaluation metrics $\mathrm{MAPE}$, $\mathrm{RMSE}$, and $\mathrm{MAE}$.

\begin{figure}[H]
    \centering
    \begin{subfigure}[b]{1\textwidth}
        \includegraphics[width=1\linewidth]{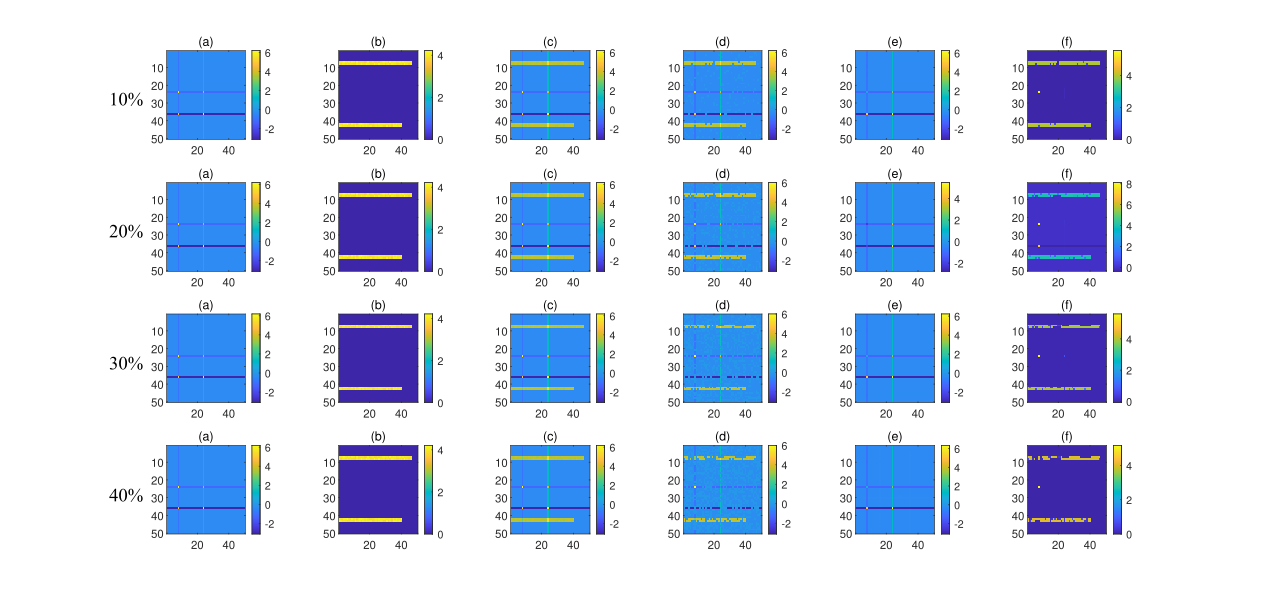} 
    \end{subfigure}
    \vspace{-1cm}
    \begin{subfigure}[b]{1\textwidth}
        \includegraphics[width=1\linewidth]{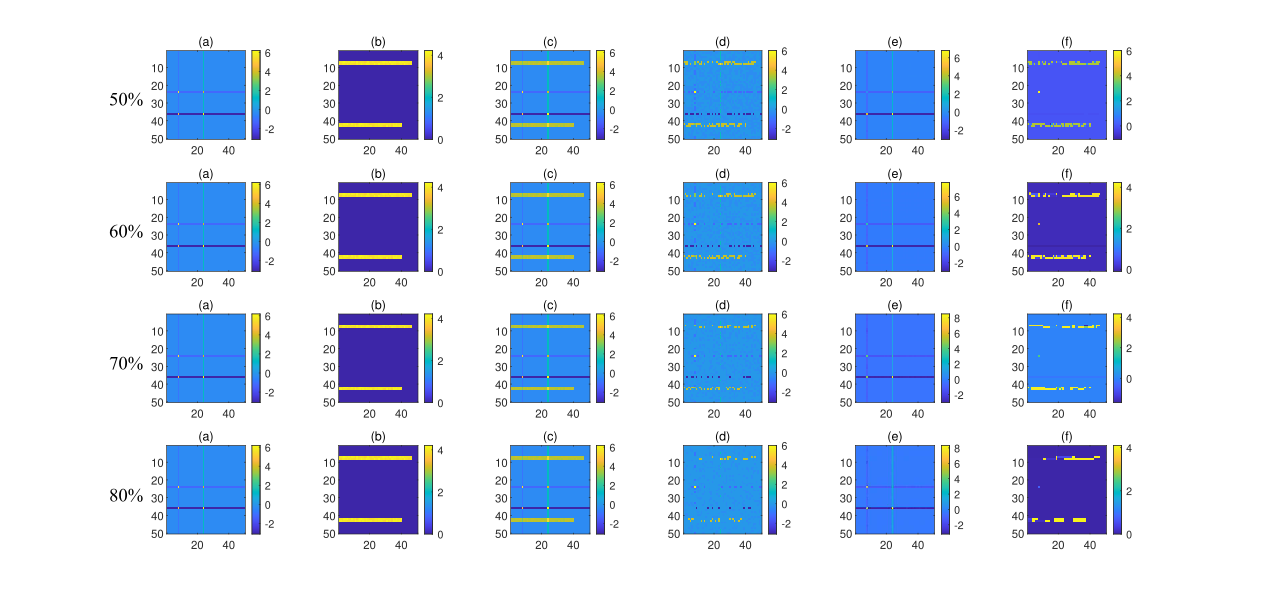} 
    \end{subfigure}
    
    \captionsetup{font=footnotesize}
    \caption{Data imputation and anomaly diagnosis performance under varying missing rates. Each row demonstrates results at 10\%, 20\%, 30\%, and 40\% missing rates respectively. Subfigure (a) shows the low-rank component of the original tensor slice, (b) displays the sparse anomaly component, (c) presents the complete original slice, (d) depicts the observed incomplete slice, (e) visualizes the recovered low-rank slice, and (f) identifies the detected block-sparse anomaly slice. All slices represent the 12th frontal slice along the third mode of the tensor. }
    \label{yichang}
\end{figure}

The Real-world dataset inherently possesses noise and anomalies, with their precise locations remaining unknown. Thus, to evaluate our model's performance on real-world dataset, we implement smoothing processing on the dataset. The specific procedure involves: (i) We select an appropriate Tucker rank to decompose the dataset, where [2,5,6] is chosen as the Tucker rank based on the criterion\cite{luo}: 
\begin{equation}r_{n}=\min_{1\leqslant l\leqslant l_{n}}\left\{l:\frac{\sum_{i=1}^{l}\tilde{\sigma}_{i}^{2}\left(\mathbf{X}_{(n)}\right)}{\sum_{i=1}^{l_{n}}\tilde{\sigma}_{i}^{2}\left(\mathbf{X}_{(n)}\right)}\geqslant\theta\right\},n=1,2,3,\end{equation}
where $\theta$ is typically set to 0.95. (ii) The core tensor and factor matrices obtained from the Tucker decomposition are utilized to reconstruct a new tensor. This reconstructed tensor preserves main features of  original data while maintaining its low-rank characteristics, thereby serving as the ground truth for our numerical experiments.

Anomalies and missing entries are introduced into the smoothed data to generate the observed tensor. Before feeding the observed data into the model, we apply an invertible scaling transformation to the observed tensor to enhance its compatibility with the model parameters. The specific procedure involves: (i) The observed tensor is element-wise subtracted by 20 and then decomposed via Tucker decomposition with rank [2,5,6]. Then we get a core tensor and factor matrices. (ii) The core tensor is element-wise multiplied by 0.14 and then reconstructed with the factor matrices to form a new tensor, which serves as the observed tensor input for the model.

\begin{table}[H]
\centering
\caption{Anomaly Diagnosis and Data Completion Evaluation on GUANGZHOU Dataset}
\begin{tabular}{ccccccccc}
\toprule
Missing rate          & 0.1      & 0.2      & 0.3      & 0.4      & 0.5      & 0.6      & 0.7      & 0.8           \\
\midrule
Precision & 0.987 & 0.989 & 0.960 & 0.981 & 0.963 & 0.972 & 0.979 & 0.976\\
Recall & 0.969 & 0.942 & 0.902 & 0.843 & 
0.767 & 0.791 & 0.765 & 0.713 \\
F1 Score & 0.978& 0.965 & 0.930 & 0.907 & 0.854 & 0.872 & 0.859 & 0.824 \\
\midrule
RMSE        & 0.002    & 0.013   & 0.013   & 0.010   & 0.043   & 0.021   & 0.013   & 0.020  \\
MAPE        & 0.104\%   & 0.517\%   & 0.519\%   & 0.594\%   & 1.521\%   & 1.234\%   & 1.11\%   & 2.00\%    \\
MAE & 0.0007& 0.0005 & 0.0027 & 0.0034 & 0.0089 & 0.0064 & 0.0061 & 0.0097\\
\midrule
\text{\footnotesize{Reverted data:}} & & & & & & & & \\
RMSE & 0.015 & 0.006 & 0.089 & 0.071 & 0.311 & 0.152 & 0.095 & 0.142 \\
MAPE & 0.025\% & 0.042\% & 0.097\% & 0.138\% & 0.338\% & 0.298\% & 0.291\% & 0.58\% \\
MAE & 0.0049 & 0.0035 & 0.0190 & 0.0244 & 0.0680 & 0.0460 & 0.0482 & 0.0726 \\
\bottomrule
\end{tabular}
\end{table}

\subsubsection{Ablation Study}
To validate the effectiveness of the proposed regularization terms, we conduct ablation studies to systematically evaluate their individual contributions to the overall model performance.
We sequentially remove individual terms, pairwise combinations and all of them. Then we compare their performance with our model under different missing rates on synthetic dataset.

According to Table \ref{ablation}, we observe that all three regularization terms can effectively enhance the model's capabilities in anomaly detection and tensor completion. It is noteworthy that this table also reveals an essential trade-off between completion efficacy and anomaly detection performance. The first regularization term, designed to enforce data continuity across three dimensions, demonstrates that its removal would degrade the completion performance while conversely enhancing anomaly detection outcomes. The second and third regularization terms, which respectively characterize the sporadic nature and local continuity of anomalous events, exhibit an inverse pattern: their elimination reduces anomaly detection effectiveness but improves completion performance. Compared to models with partial regularization term removals, the complete model demonstrates better balanced capability in coordinating data completion with anomaly detection tasks.
\begin{table}[hpt]
\centering
\caption{Experimental Results of Ablation Study}
\label{ablation}
\begin{threeparttable}
\begin{tabular}{@{}cccccccccc@{}}
\toprule
\multicolumn{2}{c}{Missing rate}                            & 0.1       & 0.2       & 0.3       & 0.4       & 0.5       & 0.6       & 0.7       & 0.8       \\
\midrule
\multicolumn{1}{c|}{\multirow{8}{*}{F1 Score}}  & a\tnote{1}        & 0.918   & \textbf{0.805}     & \textbf{0.921}     & \textbf{0.869}     & 0.829     & 0.776     & 0.751     & 0.713     \\
\multicolumn{1}{c|}{}                           &b\tnote{2}         & 0.178     & 0.178     & 0.178     & 0.178     & 0.178     & 0.178     & 0.178     & 0.178     \\
\multicolumn{1}{c|}{}                           & c\tnote{3}          & 0.728     & 0.633     & 0.683     & 0.673     & 0.600     & 0.517     & 0.399     & 0.314     \\
\multicolumn{1}{c|}{}                           & d\tnote{4}          & 0.178     & 0.178     & 0.178     & 0.178     & 0.178     & 0.178     & 0.178     & 0.178     \\
\multicolumn{1}{c|}{}                           & e\tnote{5}          & 0.728     & 0.667     & 0.673     & 0.675     & 0.602     & 0.527     & 0.399     & 0.315     \\
\multicolumn{1}{c|}{}                           & f\tnote{6}          & 0.178     & 0.178     & 0.178     & 0.178     & 0.178     & 0.178     & 0.178     & 0.178     \\
\multicolumn{1}{c|}{}                           & g\tnote{7}          & 0.178     & 0.178     & 0.178     & 0.178     & 0.178     & 0.178     & 0.178     & 0.178     \\
\multicolumn{1}{c|}{}                           & TSLTO  &  \bf{0.936}   & \textbf{0.805}   & 0.892   & 0.858   & \textbf{0.853}   & \textbf{0.782}   & \textbf{0.755}   & \bf{0.729}         \\
\midrule
\multicolumn{1}{c|}{\multirow{8}{*}{MAPE}}      &  a        & 6.494\%   & 4.439\%   & 6.660\%   & 6.789\%   & 7.008\%   & 9.244\%   & 11.025\%  & 14.395\%  \\
\multicolumn{1}{c|}{}                           & b         & 103.539\% & 97.491\%  & 96.671\%  & 93.420\%  & 90.158\%  & 90.836\%  & 93.536\%  & 96.258\%  \\
\multicolumn{1}{c|}{}                           & c        & 6.008\%   & 9.677\%   & \textbf{4.030\%}   & 2.591\%   & 2.624\%   & 5.336\%   & 5.448\%   & 3.197\%   \\
\multicolumn{1}{c|}{}                           & d         & 103.524\% & 99.380\%  & 98.022\%  & 92.947\%  & 90.802\%  & 93.008\%  & 94.704\%  & 93.524\%  \\
\multicolumn{1}{c|}{}                           & e        & \textbf{5.853\%}   & 6.465\%   & 4.687\%   & \textbf{2.387\%}   & \textbf{2.363\%}   & \textbf{2.244\%}   & \textbf{5.295\%}   & \textbf{2.624\%}   \\
\multicolumn{1}{c|}{}                           & f       & 108.794\% & 107.684\% & 104.415\% & 102.583\% & 101.633\% & 101.466\% & 102.185\% & 101.494\% \\
\multicolumn{1}{c|}{}                           & g      & 108.543\% & 119.158\% & 105.744\% & 102.503\% & 101.595\% & 105.206\% & 101.990\% & 102.573\% \\
\multicolumn{1}{c|}{}                           & TSLTO &  5.924\% & \textbf{3.941\%} & 6.506\% & 6.168\% & 6.923\% & 9.328\% & 8.220\% & 11.229\%           \\
\midrule
\multicolumn{1}{c|}{\multirow{8}{*}{RMSE}}      & a        & 0.431     & \textbf{0.388}     & 0.416     & 0.514     & \textbf{0.427}     & \textbf{0.411}     & 0.503     & 0.505     \\
\multicolumn{1}{c|}{}                           & b    & 0.603     & 0.584     & 0.577     & 0.574     & 0.565     & 0.571     & 0.614     & 0.565     \\
\multicolumn{1}{c|}{}                           & c        & 0.612     & 0.638     & 0.603     & 0.547     & 0.517     & 0.532     & 0.549     & 0.502     \\
\multicolumn{1}{c|}{}                           & d      & 0.602     & 0.584     & 0.587     & 0.573     & 0.566     & 0.572     & 0.611     & 0.554     \\
\multicolumn{1}{c|}{}                           & e       & 0.611     & 0.604     & 0.620     & 0.546     & 0.513     & 0.503     & 0.546     & 0.481     \\
\multicolumn{1}{c|}{}                           & f      & 0.743     & 0.742     & 0.740     & 0.732     & 0.717     & 0.717     & 0.726     & 0.720     \\
\multicolumn{1}{c|}{}                           & g        & 0.742     & 0.761     & 0.761     & 0.731     & 0.716     & 0.707     & 0.724     & 0.705     \\
\multicolumn{1}{c|}{}                           & TSLTO &         \textbf{0.428}   & 0.445   & \textbf{0.392}   & \textbf{0.407}   & 0.448   & 0.455   & \textbf{0.456}   & \textbf{0.470}        \\
\midrule
\multicolumn{1}{c|}{\multirow{8}{*}{MAE}}       & a        & 0.038     & \textbf{0.036}     & 0.038     & 0.047     & 0.040     & \textbf{0.048}     & 0.058     & 0.074     \\
\multicolumn{1}{c|}{}                           & b         & 0.402     & 0.394     & 0.377     & 0.364     & 0.354     & 0.353     & 0.361     & 0.360     \\
\multicolumn{1}{c|}{}                           & c         & 0.120     & 0.147     & 0.103     & 0.079     & 0.074     & 0.078     & 0.093     & 0.067     \\
\multicolumn{1}{c|}{}                           & d        & 0.402     & 0.394     & 0.383     & 0.362     & 0.357     & 0.359     & 0.362     & 0.357     \\
\multicolumn{1}{c|}{}                           &e         & 0.120     & 0.126     & 0.114     & 0.078     & 0.071     & 0.069     & 0.092     & 0.064     \\
\multicolumn{1}{c|}{}                           & f        & 0.483     & 0.480     & 0.471     & 0.464     & 0.457     & 0.455     & 0.461     & 0.456     \\
\multicolumn{1}{c|}{}                           & g        & 0.482     & 0.509     & 0.481     & 0.463     & 0.457     & 0.460     & 0.460     & 0.453     \\
\multicolumn{1}{c|}{}                           & TSLTO &  \textbf{0.035}   & 0.039   & \textbf{0.037}   & \textbf{0.037}   & \textbf{0.039}   & 0.049   & \textbf{0.047}   & \textbf{0.060}   \\
\bottomrule
\end{tabular}
\begin{tablenotes}
\footnotesize            
    \item[1] Model without $\ell_{2,0}$-norm regulations.
    \item[2] Model without $\ell_{0}$-norm regulation of $\mathcal{R}$.
    \item[3] Model without $\ell_{0}$-norm regulation of $\boldsymbol{T}_{l}\mathcal{R}_{[1]}\boldsymbol{T}_{r}^{\top}$.
    \item[4] Model without $\ell_{2,0}$-norm regulations and $\ell_{0}$-norm regulation of $\mathcal{R}$.
    \item[5] Model without $\ell_{2,0}$-norm regulations and $\ell_{0}$-norm regulation of $\boldsymbol{T}_{l}\mathcal{R}_{[1]}\boldsymbol{T}_{r}^{\top}$.
    \item[6] Model without $\ell_{0}$-norm regulation of $\mathcal{R}$ and $\ell_{0}$-norm regulation of $\boldsymbol{T}_{l}\mathcal{R}_{[1]}\boldsymbol{T}_{r}^{\top}$.
    \item[7] Model without any regulations.
\end{tablenotes}

\end{threeparttable}
\end{table}
\subsubsection{Sensitivity Analysis}
\begin{figure}[H]
    \centering
    \begin{subfigure}[b]{0.3\textwidth}
        \includegraphics[width=1\linewidth,height=4cm]{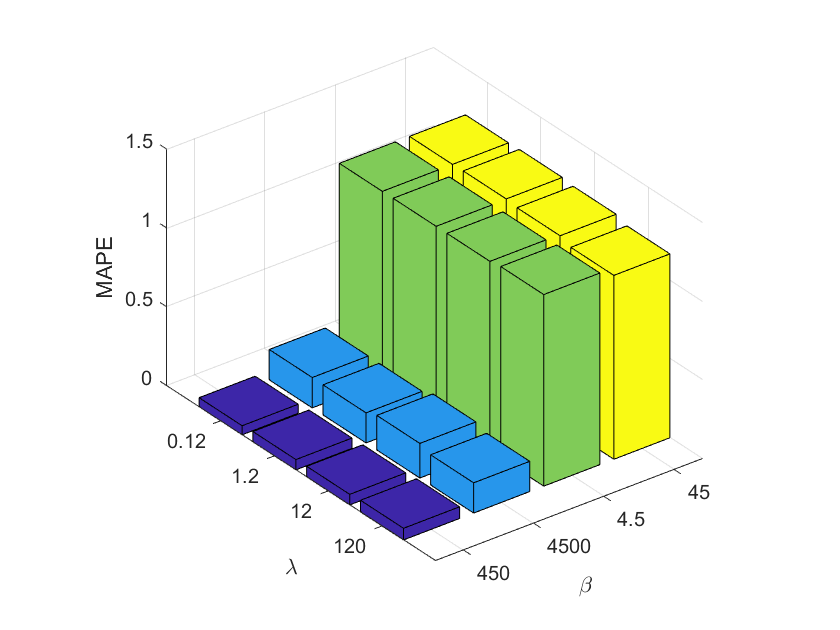}
    \caption{MAPE}
    \label{fig:enter-label}
    \end{subfigure}
    \hfill
    \begin{subfigure}[b]{0.3\textwidth}
        \includegraphics[width=1\linewidth,height=4cm]{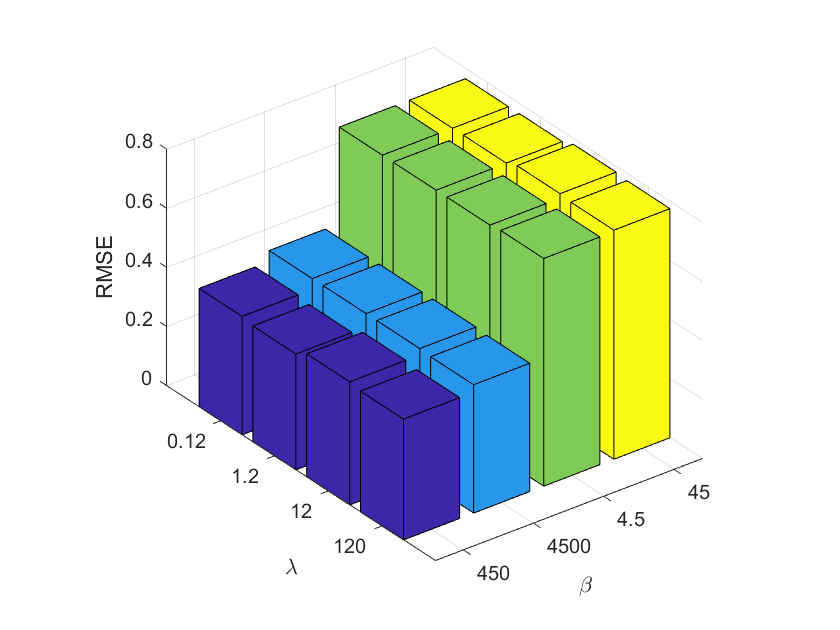}
    \caption{RMSE}
    \label{fig:enter-label}
    \end{subfigure}
    \hfill
    \begin{subfigure}[b]{0.3\textwidth}
        \includegraphics[width=1\linewidth,height=4cm]{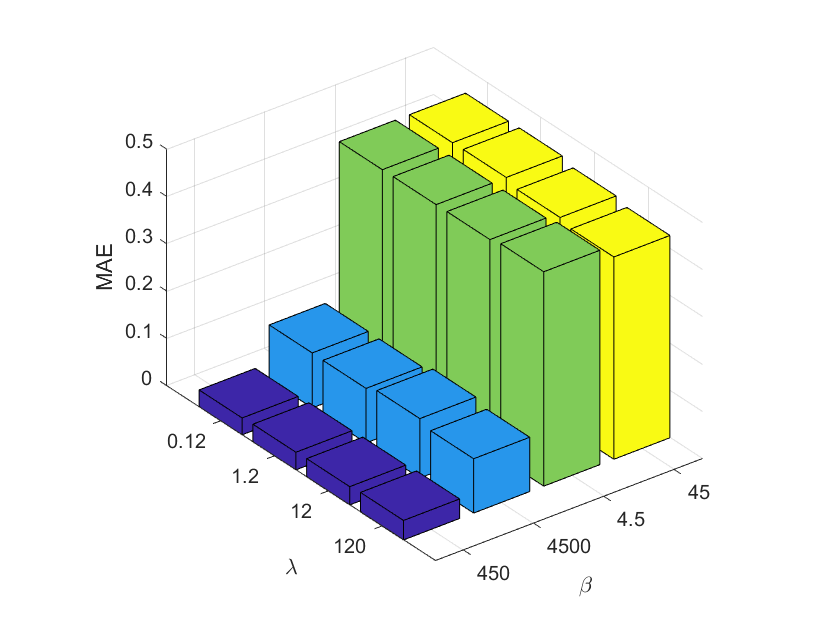}
    \caption{MAE}
    \label{fig:enter-label}
    \end{subfigure}
    \vspace{0.1cm}
    \begin{subfigure}[b]{0.3\textwidth}
        \includegraphics[width=1\linewidth,height=4cm]{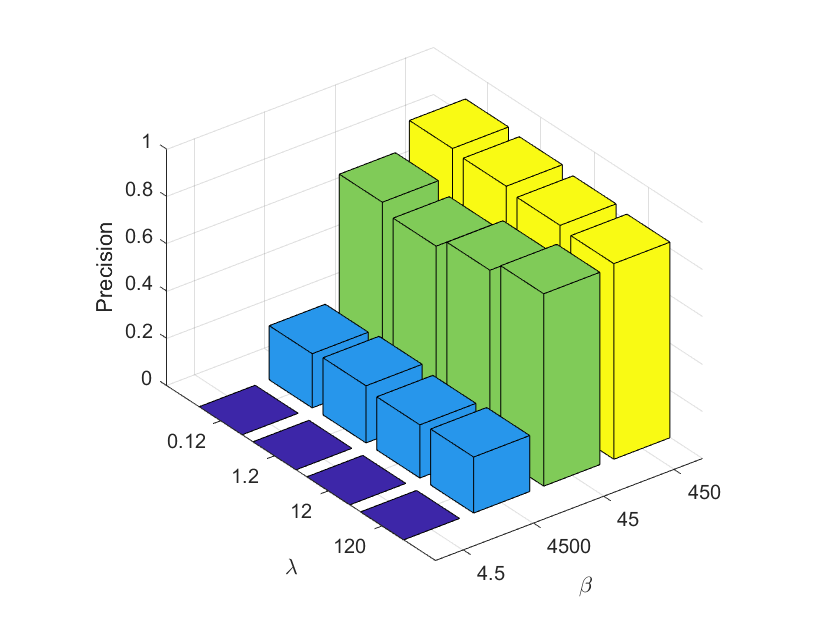}
    \caption{Precision}
    \label{fig:enter-label}
    \end{subfigure}
    \hfill
    \begin{subfigure}[b]{0.3\textwidth}
        \includegraphics[width=1\linewidth,height=4cm]{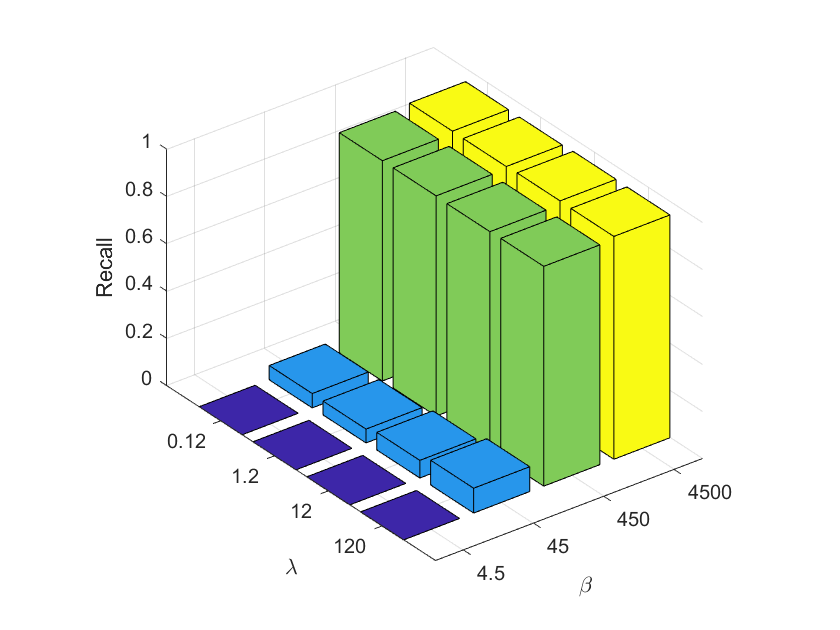}
    \caption{Recall}
    \label{fig:enter-label}
    \end{subfigure}
    \hfill
    \begin{subfigure}[b]{0.3\textwidth}
        \includegraphics[width=1\linewidth,height=4cm]{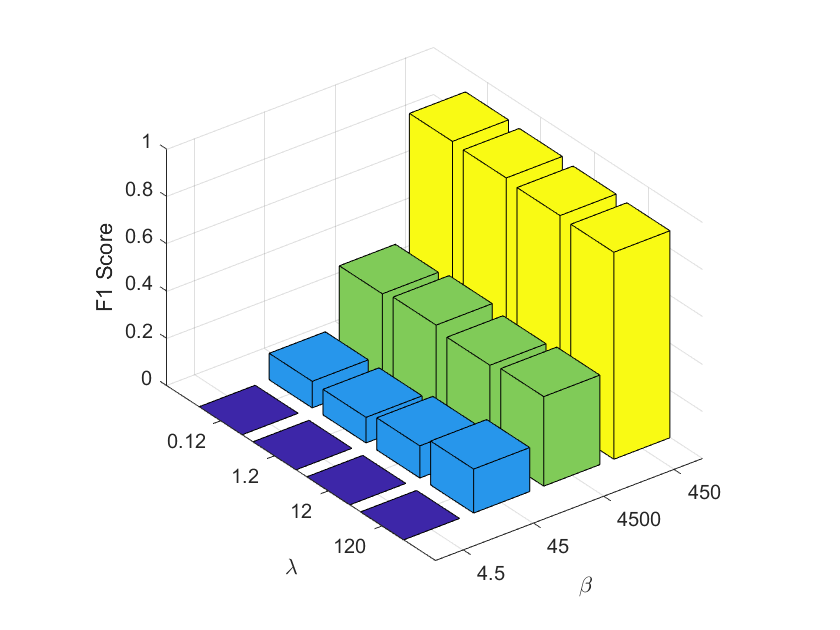}
    \caption{F1 Score}
    \label{fig:enter-label}
    \end{subfigure}
    \caption{Grid Search Results of $\beta$ and $\lambda$}
    \label{nya}
\end{figure}

\begin{figure}[H]
    \centering
    \begin{subfigure}[b]{0.3\textwidth}
        \includegraphics[width=1\linewidth,height=3.5cm]{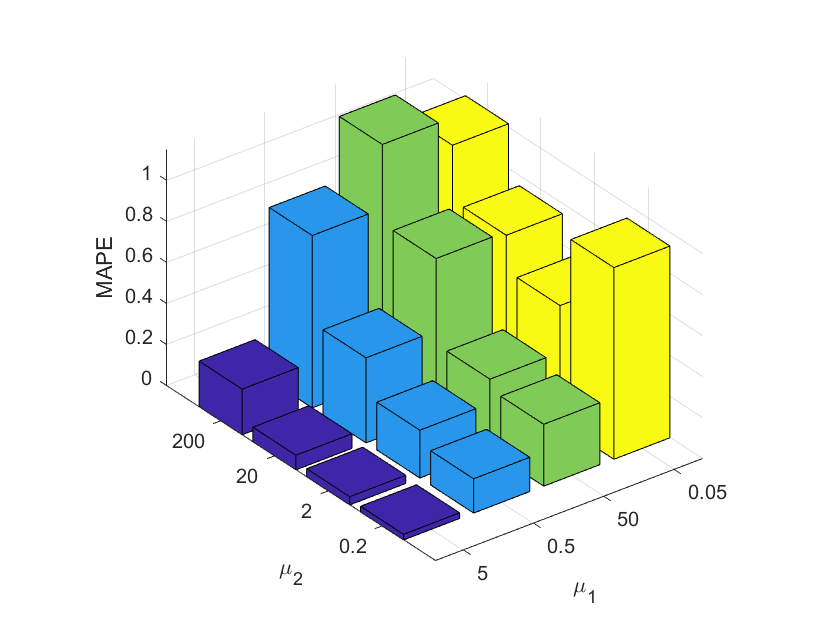}
    \caption{MAPE}
    \label{fig:enter-label}
    \end{subfigure}
    \hfill
    \begin{subfigure}[b]{0.3\textwidth}
        \includegraphics[width=1\linewidth,height=3.5cm]{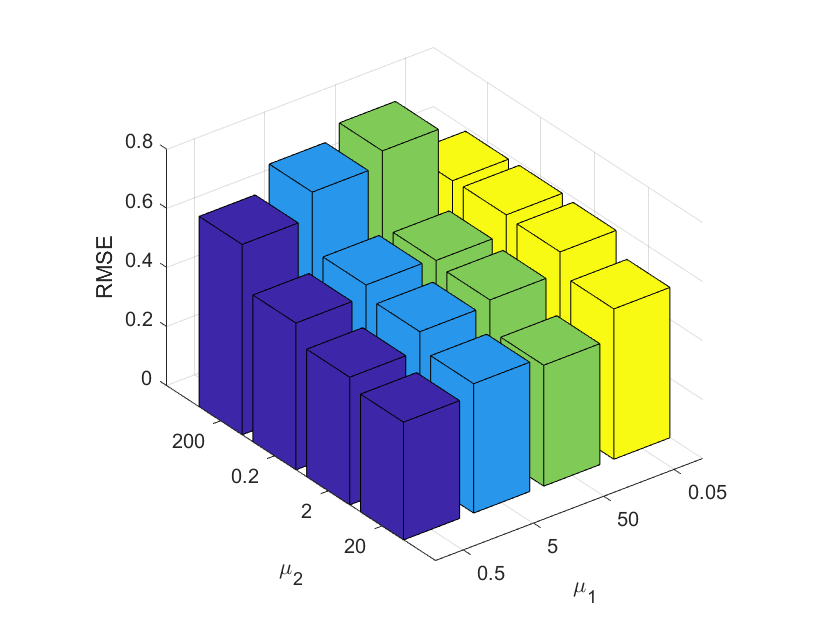}
    \caption{RMSE}
    \label{fig:enter-label}
    \end{subfigure}
    \hfill
    \begin{subfigure}[b]{0.3\textwidth}
        \includegraphics[width=1\linewidth,height=3.5cm]{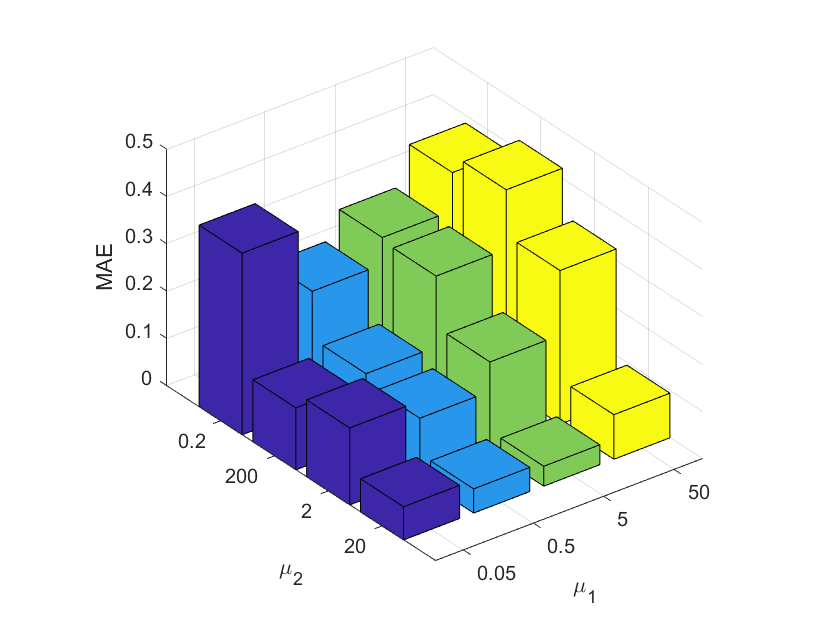}
    \caption{MAE}
    \label{fig:enter-label}
    \end{subfigure}
    \vspace{0.1cm}
    \begin{subfigure}[b]{0.3\textwidth}
        \includegraphics[width=1\linewidth,height=3.5cm]{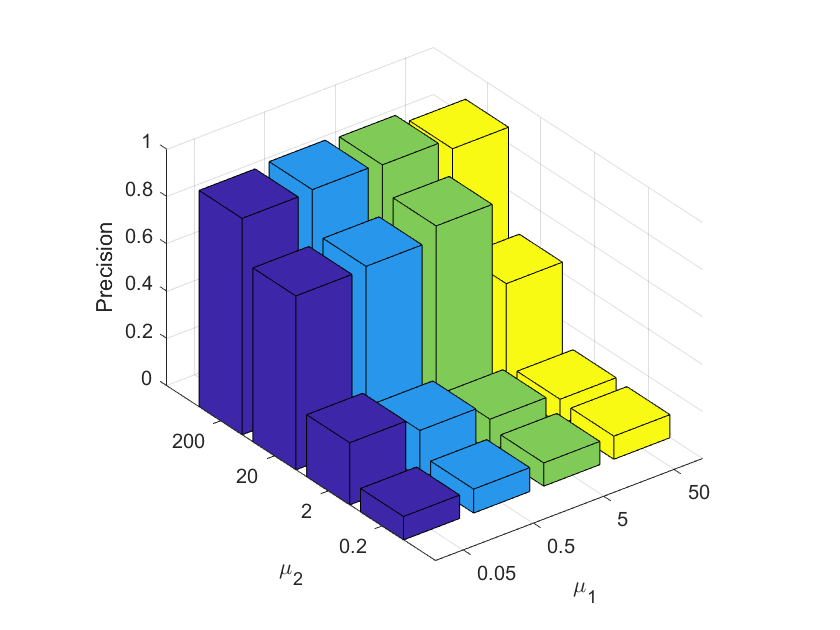}
    \caption{Precision}
    \label{fig:enter-label}
    \end{subfigure}
    \hfill
    \begin{subfigure}[b]{0.3\textwidth}
        \includegraphics[width=1\linewidth,height=3.5cm]{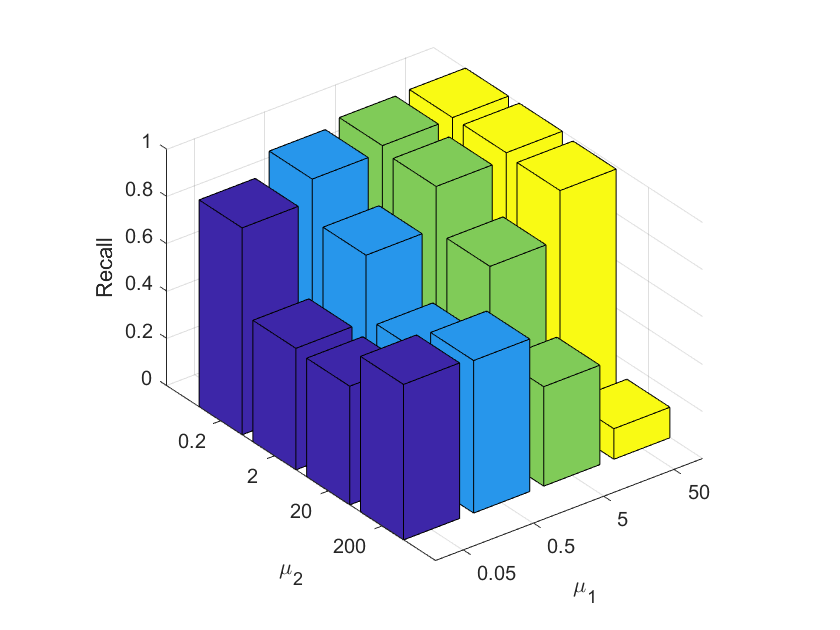}
    \caption{Recall}
    \label{fig:enter-label}
    \end{subfigure}
    \hfill
    \begin{subfigure}[b]{0.3\textwidth}
        \includegraphics[width=1\linewidth,height=3.5cm]{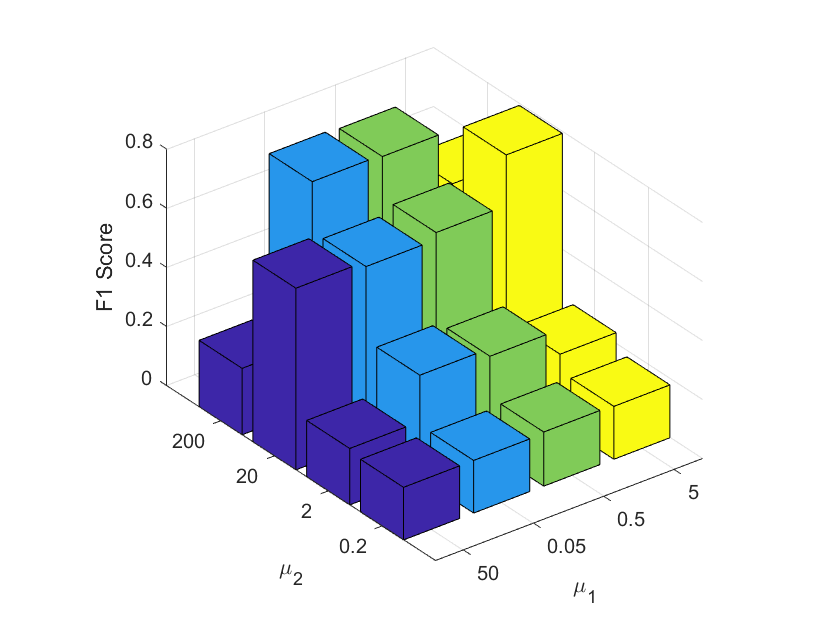}
    \caption{F1 Score}
    \label{fig:enter-label}
    \end{subfigure}
    \caption{Grid Search Results of $\mu_1$ and $\mu_2$}
    \label{miao}
\end{figure}

In 4.3.1, we demonstrated that all regularization terms contribute to enhancing model performance. We further investigate the parameter interplay between regularization terms and the least squares component. Given that $\beta$ and $\lambda$ characterize tensor recovery performance while $\mu_{1}$ and $\mu_{2}$ characterize anomaly detection performance, we conduct a grid search on simulated data with 30\% missing entries for these two parameter groups respectively. This systematic exploration, as shown in Figure \ref{nya} and Figure \ref{miao}, can reveal their joint effects on both complementary tasks. 

\subsection{Algorithm Comparison}

In this section, we present a comparative analysis of the proposed model and RMC21\cite{huyue}\footnote{\url{https://github.com/Lab-Work/Robust_tensor_recovery_for_traffic_events}} (Robust Tensor Recovery with Fiber Outliers) in terms of performance for traffic flow data completion and anomaly detection tasks.
The RMC21 models the global low-rank structure through Tucker decomposition, identifies sparse anomalous fibers via the $\ell_{2,1}$ -norm regularization, and performs missing data imputation specifically on unobserved entries within non-anomalous regions.

\begin{table}[H]
\centering
\caption{Anomaly Diagnosis and Data Completion Comparison}
\label{comparison}
\begin{threeparttable}
\begin{tabular}{@{}cccccccccc@{}}
\toprule
\multicolumn{2}{c}{Missing rate}        & 0.1     & 0.2     & 0.3     & 0.4     & 0.5     & 0.6     & 0.7     & 0.8      \\
\midrule
\multicolumn{1}{c|}{\multirow{2}{*}{Precision}} & RMC21         & 0.697    & 0.697    & 0.693    & 0.687    & 0.678    & 0.652    & 0.609    & 0.576      \\
\multicolumn{1}{c|}{}  & TSLTO & 0.987 & 0.989 & 0.960 & 0.981 & 0.963 & 0.972 & 0.979 & 0.976\\
\multicolumn{1}{c|}{\multirow{2}{*}{Recall}} & RMC21          & 0.883    & 0.748    & 0.603    & 0.595    & 0.478    & 0.346    & 0.296    & 0.202        \\
\multicolumn{1}{c|}{}  & TSLTO& 0.969 & 0.942 & 0.902 & 0.843 & 
0.767 & 0.791 & 0.765 & 0.713 \\
\multicolumn{1}{c|}{\multirow{2}{*}{F1 Score}} & RMC21         & 0.779    & 0.721    & 0.645    & 0.638    & 0.560    & 0.452    & 0.398    & 0.299        \\
\multicolumn{1}{c|}{}  & TSLTO &0.978& 0.965 & 0.930 & 0.907 & 0.854 & 0.872 & 0.859 & 0.824\\
\midrule
\multicolumn{1}{c|}{\multirow{2}{*}{rmse\tnote{1}}} & RMC21        & 0.658    & 0.643    & 0.614    & 0.678    & 0.678    & 0.664    & 0.793    & 0.849        \\
\multicolumn{1}{c|}{}  & TSLTO & 0.002       & 0.012       & 0.008      & 0.070       & 0.025       & 0.013       & 0.093       & 0.014        \\
\multicolumn{1}{c|}{\multirow{2}{*}{mape\tnote{2}}} &RMC21        & 4.324\%  & 4.242\%  & 4.050\%  & 5.116\%  & 5.157\%  & 5.622\%  & 8.798\%  & 12.804\%        \\
\multicolumn{1}{c|}{}  & TSLTO & 0.076\%       & 0.509\%       & 0.400\%       & 0.472\%       & 1.111\%       & 0.960\%       & 0.910\%       & 1.661\%        \\
\multicolumn{1}{c|}{\multirow{2}{*}{mae\tnote{3}}} & RMC21       & 0.124    & 0.119    & 0.110    & 0.133    & 0.135    & 0.137    & 0.196    & 0.248        \\
\multicolumn{1}{c|}{}  & TSLTO & 0.0007       & 0.0005       & 0.0020       & 0.0026       & 0.0062       & 0.0047       & 0.0048       & 0.0079        \\
\midrule
\multicolumn{1}{c|}{\multirow{2}{*}{RMSE}} & RMC21          & 1.095    & 1.095    & 1.096    & 1.102    & 1.105    & 1.114    & 1.165    & 1.199       \\
\multicolumn{1}{c|}{}  & TSLTO & 0.002    & 0.013   & 0.013   & 0.010   & 0.043   & 0.021   & 0.013   & 0.020 \\
\multicolumn{1}{c|}{\multirow{2}{*}{MAPE}} & RMC21          & 13.449\% & 14.697\% & 17.282\% & 14.202\% & 14.329\% & 16.378\% & 17.430\% & 21.084\%        \\
\multicolumn{1}{c|}{}  & TSLTO& 0.104\%   & 0.517\%   & 0.519\%   & 0.594\%   & 1.521\%   & 1.234\%   & 1.11\%   & 2.00\%\\
\multicolumn{1}{c|}{\multirow{2}{*}{MAE}} & RMC21          & 0.369    & 0.367    & 0.364    & 0.375    & 0.380    & 0.390    & 0.433    & 0.479       \\
\multicolumn{1}{c|}{}  & TSLTO& 0.0007& 0.0005 & 0.0027 & 0.0034 & 0.0089 & 0.0064 & 0.0061 & 0.0097\\
\bottomrule
\end{tabular}
\begin{tablenotes}
    \item[1] rmse only meassures RMSE at non-anomalous entries.
    \item[2] mape only meassures MAPE at non-anomalous entries.
    \item[3] mae only meassures MAE at non-anomalous entries.
\end{tablenotes}
\end{threeparttable}
\end{table}

\section{CONCLUSION}
A sparse and low-rank tensor optimization model has been proposed for traffic flow data recovery and anomaly diagnosis, where the element-wise sparsity, row sparsity and Tucker low-rankness have been adopted to effectively exploit the sparsity of anomalies and spatiotemporal correlations in traffic flow data. An efficient ADMM approach has been designed to handle the proposed nonconvex discontinuous optimization model, and numerical experiments have been demonstrated the effectivity of our proposed approach. In the future research, the traffic speed prediction, along with the congestion analysis, will be further considered, and some distributed second-order stochastic optimization methods will be studied to accelerate the computation for large-scale intelligent traffic systems.

\bibliographystyle{apalike}  
\bibliography{main}  

\end{document}